\documentclass[runningheads]{llncs}
\usepackage[utf8]{inputenc}

\usepackage[pdftex,colorlinks,citecolor=blue,urlcolor=blue]{hyperref}
\usepackage{amssymb}
\usepackage{blkarray}
\usepackage{graphicx,color,tikz,pgf}
\usepackage{amsmath}
\usepackage[ruled]{algorithm2e}
\usepackage{comment}
\usepackage{cleveref}
\newcommand{\init}{\mathrm{in}}
\graphicspath{{./Figures_Indep_Model/}}
\usepackage{array}

\begin{document}

\title{Optimal Transport to a Variety}

\author{T\"urk\"u \"Ozl\"um \c{C}elik \inst{1}
\and
Asgar Jamneshan \inst{3}
\and
Guido Mont\'ufar \inst{1,3} 
\and \\
Bernd Sturmfels \inst{1,2}
\and
Lorenzo Venturello \inst{1}
}
\authorrunning{T. \"O. \c{C}elik, A. Jamneshan, G. Mont\'ufar, B. Sturmfels, L. Venturello }
\institute{Max Planck Institute for Mathematics in the Sciences, Leipzig \and
University of California at Berkeley \and University of California at Los Angeles
}

\maketitle    

\begin{abstract}
We study the problem of minimizing the Wasserstein distance between a probability distribution and an algebraic variety. 
We consider the setting of finite state spaces and describe the solution depending on the choice of the ground metric and the given distribution. 
The Wasserstein distance between the distribution and the variety is the minimum of a linear functional over a union of transportation polytopes. 
We obtain a description in terms of the solutions of a finite number of
systems of polynomial equations. 
The case analysis is based on the ground metric. 
A detailed analysis is given for the two bit independence~model.

\keywords{Algebraic Statistics \and Linear Programming \and 
Optimal Transport Estimator \and Polynomial Optimization  \and
Transportation Polytope \and Triangulation \and Wasserstein Distance}
\end{abstract}

\section{Introduction} 
Density estimation in statistics is the problem of learning a hypothesis density $\nu$ 
based on samples $x_1,\ldots,x_N\in\Omega$ from an unknown density $\mu$. A standard approach to solving this problem is~to define a statistical model
$\mathcal{M}$ of candidate hypotheses, and then select a density
from $\mathcal{M}$ that minimizes some type of distance to the empirical distribution
 $\bar \mu = \frac{1}{N}\sum_i \delta_{x_i}$. 
An example of this is the maximum likelihood estimator
\cite[Chapter 7]{sullivant2018algebraic}, which minimizes the Kullback-Leibler divergence between $\bar\mu$ and $\mathcal{M}$.
 This estimator selects $\nu\in\mathcal{M}$ by maximizing the
 log-likelihood  $ \sum_{i=1}^N \log \nu(x_i)$. 

When the sample space $\Omega$ is a metric space, optimal transport defines a distance between probability distributions \cite{villani08}. The corresponding estimator selects $\nu\in\mathcal{M}$ so that it assigns a high probability to points $x$ that are close, but not necessarily equal, to samples $x_i$. 
In contrast to the maximum likelihood estimator, this incorporates 
the metric on $\Omega$.
 One key advantage of this is that distances between distributions are well defined  even when they have disjoint supports. 
The minimum Wasserstein distance estimator plays an important role in machine learning applications. 

The key disadvantage of the optimal transport distance is that it is defined as the solution to an optimization problem. Thus, computing the minimum Wasserstein distance estimator requires solving a double minimization problem. 
In a few special cases, the Wasserstein distance can be given by a formula, e.g.~ in the case of two Gaussian distributions. 
However, for general ground distances and distributions, a closed formula is not available. The standard methods for numerical computation of the Wasserstein distance between two distributions have super cubic complexity in the size of the distributions \cite{5459199}. 
Therefore, much work has been devoted to developing 
fast methods for optimal transport~\cite{MAL-073}. 
An important advance has been the introduction of entropy regularized optimal transport and iterative computations with a Sinkhorn algorithm \cite{cuturi13}, which allows for a cheaper computation and has increased the applicability of
optimal~transport. 

In large scale problems, the exact Wasserstein distance and the minimum distance estimator remain out of reach. A very successful and popular model for obtaining implicit generative models is the Wasserstein generative adversarial network \cite{arjovsky2017wasserstein}. 
This is based on the Kantorovich dual formulation of the Wasserstein-$1$ distance, as a difference of expectation values of an optimal discriminative function. 
Training (i.e.~ fitting the parameters of the model) is based on estimating the expectations by sample averages, approximating the discriminator by a neural network, and following the negative gradient of the estimated distance with respect to the model parameters. 

A number of works address the statistical complexity of estimating the optimal transport cost. 
The asymptotic behavior of the minimum Wasserstein distance estimator was studied in \cite{BASSETTI20061298} and \cite{e50e1a8b234c4106af1f9dd6ef8fea08}. 
The convergence of the empirical distribution for increasing sample size was studied in
\cite{weed17}. 

Specifying a model beforehand allows us to focus the search for a hypothesis, reducing statistical and computational complexity. 
In many cases the model is given in terms of a parametrization with a small number of parameters, thus providing a compact representation of hypotheses. It can also be specified in terms of properties of interest, such as conditional independence relations. 
This view is taken in algebraic statistics \cite{sullivant2018algebraic}. 
When the model is an exponential family (a toric variety), maximum likelihood estimation is a convex optimization problem. 
For some exponential families, such as decomposable hierarchical models, the maximum likelihood estimator can be written explicitly
 (e.g.~\cite[Chapter 7]{sullivant2018algebraic}). Recent work characterizes 
 such cases where the solution is rational \cite{duarte2019discrete}. Closed formulas 
 are also known for some
 latent variable graphical models  \cite{MLElatentclass,seigal2018mixtures}.  

The present study is cast on the discrete side
of algebraic statistics \cite{sullivant2018algebraic}.
In our setting, the model $\mathcal{M}$ is an algebraic variety inside a probability simplex.
We wish to understand fundamental properties of the
minimum Wasserstein distance estimator for $\mathcal{M}$.
What is the structure of the function that computes the Wasserstein distance between a given data distribution and a point in $\mathcal{M}$? How does it change depending on the ground metric that is laid on the sample space? How does it change depending on the data distribution? How does it depend on the model? Is the minimizer unique, or are there finitely many minimizers? Can we obtain a closed formula? 

The optimal transport distance between two 
points in our simplex is the solution to a linear program over a transportation polytope. 
The optimal transport distance between a distribution and
$\mathcal{M}$ is the minimum of a linear functional over an infinite union of transportation polytopes.  Our aim is to
understand the combinatorics and geometry of
this parametric linear program.

This article is organized as follows. 
Section~\ref{sec:Wasserstein_distance} recalls the definition of the Wasserstein distance.  
It also provides the relevant background in linear programming,
geometric combinatorics, and commutative algebra.
A key insight is that the given metric on $\Omega$ induces a regular triangulation
of a product of two simplices (cf.~Theorem \ref{thm:sc}),
 and this induces a mixed polyhedral subdivision of one simplex when $\mu$ is fixed.
Section~\ref{sec:algorithm} presents our algorithm
for computing the Wasserstein distance
from a distribution $\mu$ to a model $\mathcal{M}$ in the probability simplex.
The main subroutine is the optimization of linear functions over the pieces
of $\mathcal{M}$ that arise from the mixed subdivision.

We illustrate 
Algorithm \ref{algorithm}  by working out the geometry for 
the discrete ground metric on three states. This is illustrated in
Figure \ref{fig:curve}.
In Section~\ref{sec:degree} we focus on the case 
of primary interest, namely when
the model $\mathcal{M}$ is an algebraic variety.
Here  the minimum Wasserstein distance estimator is a
piecewise algebraic function. We show how each piece
can be represented by the hypersurface that is dual to $\mathcal{M}$ in the
sense of projective geometry.
In Section~\ref{sec:independence} we undertake a detailed case study.
Namely, we determine
 the minimum Wasserstein estimator of a discrete independence model. 

\section{Geometric Combinatorics of the Wasserstein Distance}
\label{sec:Wasserstein_distance}

Let $\Delta_{n-1}=\{(p_1,\ldots,p_n) \in \mathbb{R}^n_{\geq 0}\,:\,
\sum_{i=1}^np_i=1\}$
denote the simplex of probability distributions  on the set $[n] = \{1,2,\ldots,n\}$.
We fix a symmetric $n\times n$ matrix $d=(d_{ij})$ with nonnegative entries.
In our application, the pair $([n],d)$ will be a finite metric space, so
we have $d_{ii} =0 $ and 
$d_{ik} \leq d_{ij} + d_{jk}$ for all $i,j,k$. 
We identify $\Delta_{n^2-1}$ with the set of nonnegative $n {\times} n$ matrices
whose entries sum to~$1$.

Fix two distributions 
  $\mu,\nu\in\Delta_{n-1}$. The associated {\em transportation polytope} is 
\begin{equation}
\label{eq:transportationpolytope}   
   \Pi(\mu,\nu)\,\,=\,\,
 \bigl\{\,\pi \in \Delta_{n^2-1}\,: \,
 \sum_{i=1}^n\pi_{ij}=\mu_j \,\,\text{for all $j$ \ and }
 \sum_{j=1}^n \pi_{ij}=\nu_i \,\,\text{for all $i$}\,
 \bigr\} .
\end{equation} 
   Thus, $\Pi(\mu,\nu) $ is the set of
nonnegative $n \times n$-matrices with prescribed row and column sums.
 This polytope has dimension $ (n-1)^2$, provided
  $\mu,\nu \in {\rm int}(\Delta_{n-1})$, and it is simple
  if $\mu,\nu$ are generic.

 We consider the linear programming problem
 on the transportation polytope $\Pi(\mu,\nu)$ with cost matrix $d$.
  This is known as the
{\em transportation problem} for $(\mu,\nu,d)$.
The optimal value of this linear program is known as
the  Wasserstein distance between $\mu$ and $\nu$
with respect to  $d$.
Thus, the {\em Wasserstein distance} is
\begin{equation}
\label{eq:wassersteindistance1}
    W(\mu,\nu)\quad =\,\,\min_{\pi\in \Pi(\mu,\nu)}\sum_{1\leq i,j\leq n} d_{ij}\pi_{ij}.
\end{equation}

We are interested in the following parametric version of this
linear programming problem. We fix any subset $\mathcal{M}$
of the model $\Delta_{n-1}$. This set  is our statistical model.
The {\em Wasserstein distance} between a given distribution $\mu$
and the model $\mathcal{M}$ with respect to the metric $d$ is defined to be
 \begin{equation}
\label{eq:wassersteindistance2}   
    W(\mu,\mathcal{M})\quad =\quad \min_{\nu\in \mathcal{M}}\,\min_{\pi\in \Pi(\mu,\nu)}\sum_{1\leq i,j\leq n}\! d_{ij}\pi_{ij}.
\end{equation}
Computing this quantity amounts to solving a nested optimization problem.
Namely, we are minimizing the cost function $d$ over the set
$\bigcup_{\nu \in \mathcal{M}} \Pi(\mu,\nu)$.
The constraint set can be thought of
as a bundle of transportation polytopes over the model $\mathcal{M}$.
Our goal is to understand its geometry.

\smallskip

The $2n$ linear constraints that define the transportation
polytope $\Pi(\mu,\nu)$ can be written as
 $A\pi=(\mu_1,\dots,\mu_n,\nu_1,\dots,\nu_n)^T$ for a certain
 matrix $A\in\{0,1\}^{2n\times n^2}$ of rank $2n-1$.
 The columns of this matrix are the vertices of the product of the
 standard simplices $\Delta_{n-1}\times\Delta_{n-1} \subset \mathbb{R}^n \times \mathbb{R}^n$.
 
\begin{example}\label{ex1}
    Let $n{=}4$. The  polytopes $\Pi(\mu,\nu)$ 
    are $9$-dimensional for $\mu,\nu \in {\rm int}(\Delta_3)$. They
        are the fibers of the linear map $\Delta_{15} \rightarrow \Delta_3 \times \Delta_3$
given by the matrix
    \setcounter{MaxMatrixCols}{20}
    $$A\,\,=\,\,\begin{bmatrix}
    1 & 1 & 1 & 1 & 0 & 0 & 0 & 0 & 0 & 0 & 0 & 0 & 0 & 0 & 0 & 0 \\
    0 & 0 & 0 & 0 & 1 & 1 & 1 & 1 & 0 & 0 & 0 & 0 & 0 & 0 & 0 & 0\\
    0 & 0 & 0 & 0 & 0 & 0 & 0 & 0 & 1 & 1 & 1 & 1 & 0 & 0 & 0 & 0 \\
    0 & 0 & 0 & 0 & 0 & 0 & 0 & 0 & 0 & 0 & 0 & 0 & 1 & 1 & 1 & 1 \\
    1 & 0 & 0 & 0 & 1 & 0 & 0 & 0 & 1 & 0 & 0 & 0 & 1 & 0 & 0 & 0\\
    0 & 1 & 0 & 0 & 0 & 1 & 0 & 0 & 0 & 1 & 0 & 0 & 0 & 1 & 0 & 0\\
    0 & 0 & 1 & 0 & 0 & 0 & 1 & 0 & 0 & 0 & 1 & 0 & 0 & 0 & 1 & 0\\
    0 & 0 & 0 & 1 & 0 & 0 & 0 & 1 & 0 & 0 & 0 & 1 & 0 & 0 & 0 & 1\\
    \end{bmatrix}.$$
\end{example}

Fix a generic matrix  $d \in \mathbb{R}^{n^2}$.
The optimal bases of our linear program (\ref{eq:wassersteindistance1}),
as the distributions $\mu,\nu$ range over the simplex $\Delta_{n-1} $, are the maximal 
simplices $\sigma$ in a  triangulation $\Sigma_d$ of 
the $(2n-2)$-dimensional polytope $\Delta_{n-1} \times \Delta_{n-1}$.
Combinatorially, such a basis $\sigma$ consists of the edges in a spanning
tree of the complete bipartite graph on $[n] \times [n]$.
Let $A_{\sigma}$ be the submatrix of $A$ given by the columns 
that are indexed by
$\sigma$. For $(\mu,\nu) \in \Delta_{n-1} \times \Delta_{n-1}$, there exists
a unique column vector $\pi_\sigma$   such that
$ A_\sigma \cdot \pi_\sigma    =  (\mu,\nu)^T$.
Note that the coordinates of $\pi_\sigma$ are linear functions in $(\mu,\nu)$.

Let $\tilde \pi_\sigma$ denote the matrix in $\mathbb{R}^{n^2}$
that agrees with $\pi_\sigma$ in all coordinates in $\sigma$
and is zero in all other coordinates.
Then $\tilde \pi_\sigma$ is the optimal 
vertex of $\Pi(\mu,\nu)$ for all pairs $(\mu,\nu)$ in the simplex $\sigma$.
On that $\sigma$,
 the Wasserstein distance between
our two distributions
is given by the linear function
\begin{equation}
\label{eq:linearfunction}
 (\mu,\nu) \quad \mapsto \quad
 W(\mu,\nu)\,\,\,=\sum_{1\leq i,j\leq n}  
  d_{ij}\cdot(\tilde \pi_{\sigma})_{ij}. 
  \end{equation}
  This allows us to remove the inner optimization 
  when solving (\ref{eq:wassersteindistance2}).  
   For each simplex $\sigma \in \Sigma_d$, our task is to
      minimize the linear function (\ref{eq:linearfunction}) over the intersection
  $ (\mu \times \mathcal{M} ) \cap \sigma$. Among these optimal
  solutions, one for each simplex $ \sigma \in \Sigma_d$,
  we then select the solution with the smallest optimal value.
  This is the geometric idea behind the algorithm that will be presented in the next section.

We now shift gears and we discuss the study of
triangulations of   $\Delta_{n-1}\times\Delta_{n-1}$.
This is a rich subject in geometric combinatorics, with numerous
connections to optimization, tropical geometry, 
enumerative combinatorics, representation theory, commutative algebra,
and algebraic geometry.  The triangulations which appear in our context 
are called \emph{regular triangulations}  \cite[Chapter 2]{DLRS}.
There are various different approaches for computing these objects. The one
we favor here is based on commutative algebra. Namely, we represent
our objects as initial ideals of the ideal of $2 \times 2$-minors of
an $n \times n$-matrix of unknowns. In the language of algebraic
geometry, these are the toric degenerations of the Segre
variety $\mathbb{P}^{n-1} \times \mathbb{P}^{n-1}$ in its embedding in 
the matrix space $\mathbb{P}^{n^2-1}$.

In the rest of this section we present relevant definitions and results.
We refer to \cite{MS,St_GB} for an extensive treatment of the subject. 
Fix the polynomial ring $R=K[y_{ij}: 1\leq i,j\leq n]$ over a field $K$. 
We identify nonnegative integer vectors $\alpha\in\mathbb{N}^{n\times n}$ 
with monomials $y^{\alpha}=\prod_{i=1}^n y_{ij}^{\alpha_{i,j}}$.
    Let $d\in\mathbb{R}^{n\times n}$ and $I$ an ideal in $R$.    Consider any
    polynomial $f=\sum_{\alpha\in \mathbb{N}^{n\times n}}c_{\alpha}y^{\alpha}\in I$.
    The
    {\em initial form} of $f$ is defined to be
    $\init_{d}(f) =
    \sum_{ d\cdot\alpha = {\bf d}}c_{\alpha}y^{\alpha}$ with ${\bf d}=\max\{d\cdot\alpha: c_{\alpha}\neq 0\}$ where $\cdot$ denotes the standard dot product. The 
    {\em initial ideal} of $I$ with respect to the weight matrix $d$ is the following ideal in $R$:
   $$\init_{d}(I)\,\,:=\,\,\langle\, \init_{d}(f): f\in I \,\rangle.$$
For a generic choice of $d$, this is a monomial ideal, i.e.~$\init_{d}(I)$ can be generated by 
monomials. In this case, we can compute
(in a computer algebra system) a corresponding reduced Gr\"obner basis
$\{g_1,g_2,\ldots,g_r\}$ of $I$. The initial monomials
${\rm in}_d(g_1),\ldots,{\rm in}_d(g_r)$ are minimal generators of ${\rm in}_d(I)$.

The connection to regular triangulations arises when $I$ is
a {\em toric ideal} $I_A$. This works for any nonnegative integer matrix $A$,
but we here restrict ourselves to the matrix $A$ whose columns
are~the vertices of $\Delta_{n-1} \times \Delta_{n-1}$, as in
Example \ref{ex1}. The toric ideal associated to $A$ is the
{\em determinantal~ideal}
$$I_A\,\,:=\,\,\langle\, y^{u^+}-y^{u^-}: \,u\in\ker(A)\,\rangle \,\,=\,\,
\langle \,\,\text{$2 \times 2$-minors of the $n \times n$ matrix $(y_{ij})$ } \rangle.
$$
The regular polyhedral subdivisions of the product of simplices
are encoded by the  initial ideals of the ideal $I_A$.

\begin{theorem}[Sturmfels' Correspondence] {\rm \cite[Theorem 9.4.5]{DLRS}} 
\label{thm:sc}
 There is a bijection between regular subdivisions of $\Delta_{n-1}\times\Delta_{n-1}$ 
 induced by $d$ and the ideals $\init_d(I_{A})$. Moreover, $\init_d(I_{A})$ is a monomial
  ideal if and only if the corresponding subdivision of $\Delta_{n-1}\times\Delta_{n-1}$ 
  is a triangulation.
\end{theorem}

Since the matrix $A$ is totally unimodular \cite[Exercise (9), page 72]{St_GB},
  all initial monomial ideals  $\init_d(I_{A})$ are squarefree \cite[Corollary 8.9]{St_GB}.
The desired triangulation $\Sigma_d$ is the simplicial complex whose
Stanley-Reisner ideal equals $\init_d(I_{A})$. This means that
the set $\mathcal{F}(\Sigma_d)$ of its maximal simplices 
in the triangulation is read off from the prime
decomposition of the squarefree monomial ideal:
\begin{equation}
\label{eq:primdec}
   \init_d(I_{A})\quad =\,\,
\bigcap_{\sigma\in \mathcal{F}(\Sigma_d)}  \!\!  \langle \,y_{ij}\,: \,{ij}\notin \sigma\,\rangle.    
\end{equation}
For a first illustration see \cite[Example 8.12]{St_GB},
where it is shown that the {\em diagonal initial ideal} of
the determinantal ideal $I_A$ corresponds to the
{\em staircase triangulation} of the polytope $\Delta_{n-1} \times \Delta_{n-1}$.

From the perspective of optimal transport,
what has been accomplished so far?
We wrote the Wasserstein distance between two distributions
locally as a linear function. This is the function in (\ref{eq:linearfunction}).
The region $\sigma$ inside $\Delta_{n-1} \times \Delta_{n-1}$
on which this formula is valid is a simplex. The set of these
simplices is the triangulation $\Sigma_d$.
The algebraic recipe (\ref{eq:primdec}) serves to compute this.
Thus, the associated primes of ${\rm in}_d(I_A)$ are the
 linear formulas for the Wasserstein distance.

\begin{remark}
If the matrix $d$ is special then
$\init_d(I_{A})$ may not be a monomial ideal.
This happens for the discrete metric on $[n]$ when $n \geq 4$.
In such a case, we break ties with a term order
to get a triangulation. Geometrically, this corresponds to replacing
$d$ by a nearby generic matrix $d_\epsilon$. However, since the optimal value function of a linear program is piecewise linear and continuous, the limit of the optimal values for $d_{\epsilon}$ as $\epsilon\rightarrow 0$ is the optimal value for $d$.
\end{remark}

The discussion above is concerned with the piecewise-linear structure of the Wasserstein
distance $W(\mu,\nu)$ when $d$ is fixed and $\mu,\nu$ vary.
The story becomes more interesting when we allow the
matrix $d$ to vary over $\mathbb{R}^{n^2}$. This brings us
to the theory of {\em secondary polytopes}.
Two generic matrices $d$ and $d'$ are considered 
{\em equivalent} if their triangulations coincide: $\Sigma_d = \Sigma_{d'}$.
The equivalence classes are open convex polyhedral cones that partition
 $\mathbb{R}^{n^2}$. This partition is the
 {\em secondary fan} of our product of simplices. This fan is the normal fan of the secondary polytope
$\Sigma(\Delta_{n-1} \times \Delta_{n-1})$, which is the Newton polytope of the product of all subdeterminants (all sizes) of the matrix $(y_{ij})$.

For a given generic matrix $d$, its equivalence class
(a.k.a.~secondary cone) can be 
read off from the reduced Gr\"obner basis
$\{g_1,g_2,\ldots,g_r\}$ of $I_A$ with respect to $d$.
The Gr\"obner basis elements are binomials $g_i = y^{u_i^+} - y^{u_i^-}$,
where $u_1,u_2,\ldots,u_r \in \mathbb{Z}^{n^2}$.
Then the desired secondary cone equals
\begin{equation}
\label{eq:secondarycone} \bigl\{ \,d \in \mathbb{R}^{n^2} \,:\,  d \cdot u_i > 0 \,\,\,
{\rm for}\,\, i=1,2,\ldots,r \bigr\}. 
\end{equation}

\begin{example} \label{ex:discrete3}
Let $n=3$ and fix the {\em discrete metric} $d \in \{0,1\}^{3\times 3}$,
which has $d_{ii} = 0$ and $d_{ij} = 1$ if $i \not=j$.
This matrix looks special but it is actually generic. 
The corresponding Gr\"obner basis equals
$$ \begin{matrix}
\{\,
\underline{y_{12} y_{21}} - y_{11} y_{22},\,
\underline{y_{12} y_{23}} - y_{13} y_{22}, \,
\underline{y_{12} y_{31}} - y_{11} y_{32},  \,
\underline{y_{13} y_{21}} - y_{11} y_{23},\,
\underline{y_{13} y_{31}} \\
- y_{11} y_{33},\,
\underline{y_{13} y_{32}} - y_{12} y_{33},  \,
\underline{y_{21} y_{32}} - y_{22} y_{31}, \,
\underline{y_{23} y_{31}} - y_{21} y_{33}, \,
\underline{y_{23} y_{32}} - y_{22} y_{33}\,\}.
\end{matrix}
$$
The initial monomials are underlined, so the secondary cone is defined by 
$$ d_{12} +  d_{21} > d_{11} +  d_{22},\,\,
d_{12} +  d_{23} > d_{13} + d_{22}, \,
\,\ldots \,,\,\,
d_{23} + d_{32} >  d_{22} + d_{33}.
$$
For any matrix in that secondary cone in $\mathbb{R}^{3 \times 3}$, 
the initial monomial ideal equals
\begin{equation}
\label{eq:primarydecomp}
\begin{matrix}
{\rm in}_d(I_A) \quad & = &\,\,\,
\langle y_{12},y_{13},y_{21},y_{23} \rangle \,\cap\,
\langle y_{12},y_{13},y_{23},y_{32} \rangle \,\cap\,
\langle y_{12},y_{13},y_{31},y_{32} \rangle \\ & & \cap\,
\langle y_{12},y_{21},y_{31},y_{32} \rangle \,\cap\,
\langle y_{13},y_{21},y_{23},y_{31} \rangle \,\cap\,
\langle y_{21},y_{23},y_{31},y_{32} \rangle .
\end{matrix}
\end{equation}
This encodes the six $4$-simplices that form the triangulation $\Sigma_d$
of  $\Delta_2 \times \Delta_2$.
\end{example}	

\section{An Algorithm and the Geometry of Triangles}
\label{sec:algorithm}
We next present our algorithm for computing the Wasserstein distance 
to a model, $W(\mu,\mathcal{M})$.
Here the model $\mathcal{M}$ is any subset of $\Delta_{n-1}$.
Our only assumption is that we have a subroutine for minimizing a linear function over
intersections of $\mu \times \mathcal{M}$ with subpolytopes $\sigma$ of $\Delta_{n-1}\times \Delta_{n-1}$.
 The case of primary interest,
when $\mathcal{M}$ is an algebraic variety, will be addressed in the next section.
We begin by giving an informal summary.
The precise version appears in Algorithm~\ref{algorithm}.

\medskip

\begin{algorithm}[H]\label{algorithm2}
    \KwIn{ 
    An $n\times n$ matrix $d=(d_{ij})$, 
    a model $\mathcal{M} \subset \Delta_{n-1}$,
    and a distribution $\mu\in\Delta_{n-1}$.}
    \KwSty{Steps 1-3:} Compute the triangulation of the polytope
    $\Delta_{n-1} \times \Delta_{n-1}$ that is given by $d$. \\
        \KwSty{Step 4:} Incorporate $\mu$ and express matrix entries
        as linear functions in $\nu \in \mathcal{M}$.\\
    \KwSty{Step 5:} For each piece, minimize a linear function over 
    the relevant part of the model~$\mathcal{M}$. \\
          \KwSty{Steps 6-7:} The smallest minimum found in Step 5
           is the {\em Wasserstein distance}~$\,W(\mu,\mathcal{M})$.
    \caption{A friendly description of the steps in Algorithm \ref{algorithm}}
\end{algorithm}

\medskip

The first step in our algorithm is the computation of the regular
triangulation $\Sigma_d$. This is done using the algebraic method
described in Section \ref{sec:Wasserstein_distance}.
As before, $I_A$ denotes the ideal of $2  \times 2$ minors
of an $n \times n$ matrix of unknowns $y = (y_{ij})$. 
The computation of $\Sigma_d$  is a preprocessing step
that depends only on $d$. Once the triangulation is known,
we can use it to treat different models $\mathcal{M}$ and different distributions $\mu$,
by starting from Step 5 of Algorithm \ref{algorithm}.

There are two sources of complexity in Algorithm \ref{algorithm}.
First, there is the subroutine in Step 5, where we minimize a linear function
over the model $\mathcal{M}$, subject to 
nonnegativity constraints that specify $(\mu \times \Delta_{n-1}) \,\cap \,\sigma$. 
When $\mathcal{M}$ is a semialgebraic set, this is a polynomial optimization problem. For an introduction to current methods see \cite{L}. In Section \ref{sec:degree} we disregard inequality constraints and focus on the case when the model $\mathcal{M}$ is a variety. Here the complexity
is governed by the algebraic degree, which refers to the
 number of complex critical points. 
The other source of complexity is combinatorial, and it is governed
by the number of maximal simplices in the triangulation of
$\Delta_{n-1} \times \Delta_{n-1}$. This number is independent
of the triangulation. We have
\begin{equation}
\label{eq:stirling}
 \qquad \qquad |\mathcal{F}(\Sigma_d)| \,\,= \,\,
\binom{2n-2 }{n-1} \quad = \quad O\left(4^{n}n^{-1/2}\right).
\end{equation}
The second equation rests on Stirling's formula. 
This exponential complexity can be reduced when we deal
with specific finite metric spaces. Namely, if
$d$ is a symmetric matrix with very special structure,
then $\Sigma_d$ will not be a triangulation but a coarser
subdivision with far fewer cells than $\binom{2n-2}{n-1}$.
This structure can be exploited systematically, in order to
gain a reduction in complexity.

 \medskip
 
\begin{algorithm}[H]\label{algorithm}
    \KwIn{ 
    An $n\times n$ matrix $d=(d_{ij})$, 
    a model $\mathcal{M} \subset \Delta_{n-1}$,
    and  $\mu\in\Delta_{n-1}$.}
    \KwOut{The Wasserstein distance  $W(\mu,\mathcal{M})$ and a point
    in $\mathcal{M}$ that attains this distance.}
    \KwSty{Step 1:} Compute the initial $\init_{d}(I_A)$ for the ideal 
    $I_A$ of $2 \times 2$-minors. \\
    \KwSty{Step 2:} If $\init_{d}(I_A)$ is not a monomial ideal, then
    redo Step 1 with a nearby generic matrix. \\
    \KwSty{Step 3}: Compute the set $\mathcal{F}(\Sigma_d)$ of maximal simplices
    in $\Sigma_d$    using  (\ref{eq:primdec}). \\
    \KwSty{Step 4:} For every $\sigma\in \mathcal{F}(\Sigma_d)$,
    compute the matrix $\tilde \pi_\sigma$ whose entries are linear in $\nu \in \mathcal{M}$. \\
    \KwSty{Step 5:} For every $\sigma\in \mathcal{F}(\Sigma_d)$, compute
    the minimum  of the  linear function
    in   (\ref{eq:linearfunction})
    over the intersection $(\mu \times \mathcal{M}) \,\cap \,\sigma$.
       \\  
          \KwSty{Step 6:} Choose the minimum value among the 
          optimal values in Step 5. \\
  \KwSty{Step 7:}  Output this value  and
  $\nu^* \in \mathcal{M}$ satisfying $ W(\mu,\nu^*) = W(\mu,\mathcal{M})$.
    \caption{Computing the Wasserstein distance to a model}
\end{algorithm}

\smallskip

\begin{example} \label{ex:discrete4}
Consider the discrete metric $d \in \{0,1\}^{n\times n}$,
which has $d_{ii} = 0$ and $d_{ij} = 1$ if $i \not=j$.
The subdivision $\Sigma_d$ of $\Delta_n \times \Delta_n$
has $2^n-2$ maximal cells. So, it is not a triangulation for $n \geq 4$.
Combinatorially, $\Sigma_d$ is dual to 
the zonotope that is obtained by taking the Minkowski sum of
$n$ line segments in $\mathbb{R}^{n-1}$. 
This follows from the identification of triangulations of
products of simplices with tropical polytopes. The 
tropical polytope representing the discrete metric is the
$(n-1)$-dimensional {\em pyrope}; see \cite[Equation (4)]{JK}.
For instance, consider the case $n=4$: the $3$-dimensional pyrope is
the rhombic dodecahedron, which has $14$ vertices, 
$24 $ edges, and $12 $ facets \cite[Figure 4]{JK}.
\end{example}

\begin{figure}[h]
\begin{center}
\vspace{-0.14in}
	\includegraphics[scale=0.64]{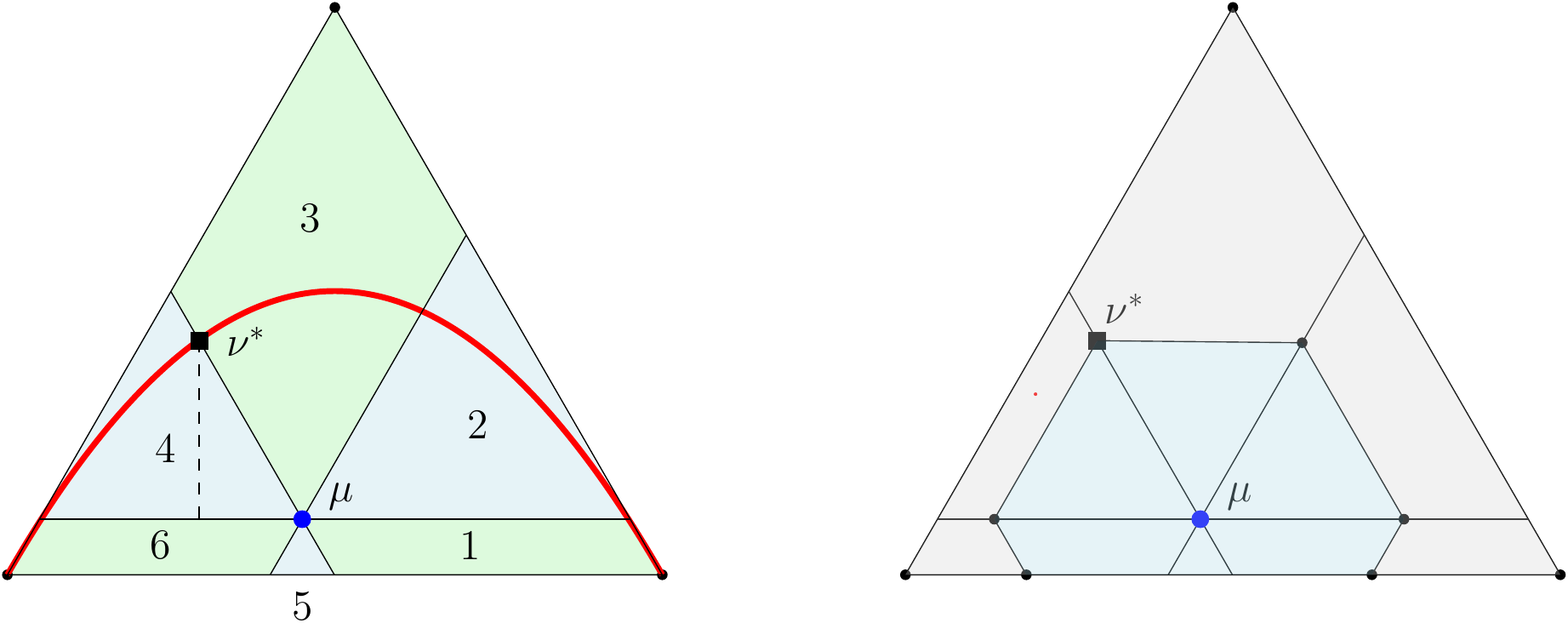}	
	\vspace{-0.1in}
	\caption{
	The model $\mathcal{M}$ is the red curve, here shown in the
	triangle $\mu \times \Delta_2$.
It intersects five of the six cells that are obtained by restricting the
	$4$-simplices $\sigma$ in $\mathcal{F}(\Sigma_d)$ 
	from $\Delta_2 \times \Delta_2$ to that triangle. The Wasserstein
	distance from $\mu$ to the curve $\mathcal{M}$ is attained by a point,
	labeled $\nu^*$, that
	lies at the intersection of two cells.
	\label{fig:curve} }
		\end{center}	
	\end{figure}

In the remainder of this section we offer a detailed illustration
of Algorithm \ref{algorithm} in the case $n=3$.
We fix the discrete metric $d$ as in Examples \ref{ex:discrete3} and \ref{ex:discrete4},
and we take $\mathcal{M}$ to be the independence model for two 
identically distributed binary random variables. This is the image of the parametrization
\begin{equation}
\label{eq:hardyweinberg}	\varphi\,:\,[0,1\,]\to\,\Delta_2\,, \quad
	p\,\mapsto \, \bigl( \,p^2,\,2p(1-p),\,(1-p)^2 \,\bigr). 
\end{equation}	
Thus $\mathcal{M} = {\rm image}(\varphi)$ is a
quadratic curve  inside the triangle $\Delta_2$.
This curve is known as the {\em Hardy-Weinberg curve} in genetics and
it is shown in red in Figure~\ref{fig:curve}.

Fix the distribution $\mu=(\,1/2,\,1/7,\,5/14\,)$.
This is marked in blue in Figure~\ref{fig:curve}. The Wasserstein
distance between $\mu$ and $\mathcal{M}$ 
is attained at $p^* = 1/\sqrt{2}$. It  equals
$$ W (\mu, \mathcal{M})\,\, =\,\, \sqrt{2} - 8/7  
\,\, = \,\,0.2713564195... \,\, = \,\,
W(\mu,\nu^*). $$
The corresponding optimal distribution $\nu^*$ in the model $\mathcal{M}$ equals
$$(\nu_1^*,\nu_2^*,\nu_3^*) \,=
\, \bigl( \,(p^*)^2,\,2p^*(1-p^*),\,(1-p^*)^2 \,\bigr) \, = \,
\bigl( \,0.5 , \,0.4142135 ..., 0.0857864... \,\bigr).$$
An optimal transportation plan is this matrix
with given row and column sums:
\begin{equation}
\label{eq:pattern}
\begin{blockarray}{cccc}
\begin{block}{(ccc)l}
\,  \pi_{11} \,& \, 0 \, & \, 0 \,\, & \; \; \nu_1^* \\
\,  \pi_{21} \,& \, \pi_{22} \, &\, \pi_{23} \,\,& \; \; \nu_2^* \\
\,   0 \,& \,0 \,&  \,\pi_{33}\, \,& \; \; \nu_3^* \\
\end{block}
\begin{block}{cccc} 
\frac{1}{2} & \frac{1}{7} & \frac{5}{14}\\
\end{block}
\end{blockarray}
 \end{equation}

This solution was found using Algorithm \ref{algorithm}.
Steps 1, 2 and 3 were already carried out in Example \ref{ex:discrete3}.
In Step 4, we translate each prime component in
 (\ref{eq:primarydecomp}) into a 
$3 \times 3$ matrix $\tilde \pi_\sigma$ whose entries are linear forms.
For instance, the third component in (\ref{eq:primarydecomp})
corresponds to the matrix in (\ref{eq:pattern}) with
$$ 
\pi_{11} =  \nu_1,\,
\pi_{21} = \mu_1 - \nu_1, \,
\pi_{22} = \mu_2 ,\,
\pi_{23} = \mu_3-\nu_3,\,
\pi_{33} = \nu_3 .
$$
We substitute
$\mu = \bigl( \frac{1}{2}, \frac{1}{7}, \frac{5}{14} \bigr)$ and
$ \nu = \bigl( \,p^2,\,2p(1-p),\,(1-p)^2 \,\bigr)$
into these six $3 \times 3$ matrices $\tilde \pi_\sigma$. 
As $\sigma$ runs over $\mathcal{F}(\Sigma_d)$, we obtain 
six feasible regions $(\mu \times \Delta_2) \cap \sigma$ in the 
$\nu$-triangle $\Delta_2$.
These are the blue triangles and the green rhombi
in Figure \ref{fig:curve}. On each of these cells,
the  objective function $\pi_{12} + \pi_{13} + \pi_{21} + \pi_{23} + \pi_{31}+ \pi_{32}$
is a quadratic function in $p$.
This quadric appears in the leftmost column of the table below,
along with the feasible region restricted to the curve $\mathcal{M}$.
The third and fourth column list the optimal solutions
that are computed in Step~5.

\begin{table}[h]
\begin{center}
\vspace{-0.14in}
  \begin{tabular}{ | c | c | c | c | }
  \hline
   Objective Function & Feasible Region & Solution $p$ & Minimum Value  \\ \hline 
  $p^2-2p+9/14$ & $0\leq p \leq \frac{(1-\sqrt{\frac{5}{7}})}{2}$   & $\frac{(1-\sqrt{\frac{5}{7}})}{2}$ &  $1/14+\sqrt{5/7}/2$  \\ \hline 
    $ -p^2+1/2$ & $\frac{(1-\sqrt{\frac{5}{7}})}{2} \leq p \leq 1-\sqrt{5/14} $ &  $1-\sqrt{5/14}$ & $2\sqrt{5/14}-6/7$  \\ \hline
  $ -2p^2+2p-\frac{1}{7}$ & $1-\sqrt{5/14}\leq p \leq \sqrt{1/2}$  & $\sqrt{1/2}$ & ${\bf \sqrt{2}-8/7}$  \\ \hline
  $-p^2+2p-9/14$ & $\sqrt{1/2}\leq p \leq \frac{(1+\sqrt{\frac{5}{7}})}{2}$ & $\sqrt{1/2}$  & ${\bf \sqrt{2}-8/7}$  \\ \hline
  $2p^2-2p+\frac{1}{7}$ & null set  & infeasible  &   \\ \hline
  $p^2-1/2$ & $\frac{(1+\sqrt{\frac{5}{7}})}{2}\leq p \leq 1$  & $\frac{(1+\sqrt{\frac{5}{7}})}{2}$ & $-1/14+\sqrt{5/7}/2$  \\
  \hline
  \end{tabular}
  \vspace{-0.2in}
  \end{center}
  \end{table}

The $i$-th row of this table corresponds to the cell labeled $i$. In Step 6 of our
algorithm, we identify cells $3$ and $4$ as those that attain the minimum value.
In Step 7 we recover the optimal solution $\nu^*$. The optimal point
is marked by $\nu^*$ in Figure \ref{fig:curve}.  The geometric fact
that the minimum is attained on the intersection of two cells
corresponds to the algebraic fact that $\pi_{21} = 0$ in~(\ref{eq:pattern}).

It is instructive to draw the balls in the Wasserstein metric
around the point $\mu$ in Figure~\ref{fig:curve}. For small radii,
these balls are regular hexagons whose sides are parallel 
to the three distinguished directions. As the radius increases, some
sides of these hexagons exit the triangle. For instance, the ball around $\mu$ that contains  optimal point 
$\nu^*$ its boundary is a non-regular hexagon containing region 5. The
boundary in each of the other regions is obtained by
drawing a line segment parallel to the opposite direction.
For instance, in region 3, we draw a horizontal segment starting
at $\nu^*$ until it hits region 2, and then we continue the boundary 
with a 60 degree turn to the right.

\section{Parametric Linear Optimization over a Variety}
\label{sec:degree}

A key step in Algorithm \ref{algorithm} is the repeated solution of linear optimization
problems over appropriate subsets of the model $\mathcal{M}$. 
We now assume that $\mathcal{M}$ is an {\em algebraic variety}
in $\Delta_{n-1} \subset \mathbb{R}^n$, i.e.~$\mathcal{M}$ consists of the nonnegative real solutions of
a system of polynomials $f_1,f_2,\dots , f_k \in \mathbb{R}[x_1,\dots , x_n]$.
We tacitly assume that $f_1 = x_1 + \cdots + x_n - 1$ is the linear equation
that cuts out the probability simplex. 
We write $X$ for the complex algebraic variety in $\mathbb{C}^n$
defined by the same equations. Let
$\bar{X}$ denote the closure of $X$ in the complex projective space~$\mathbb{P}^n$.

 When computing the Wasserstein distance from $\mu$ to the model $\mathcal{M}$ with respect to $d$,
 we must minimize a linear function over $\mathcal{M}$ subject to nonnegativity constraints
  that specify $(\mu \times \Delta_{n-1})\cap \sigma$. Here $\sigma$ runs over all
  maximal simplices in the triangulation $\Sigma_d$ of $\Delta_{n-1} \times \Delta_{n-1}$.
  Let us assume for simplicity that the minimum is attained at a smooth point of $X$
  that is in the relative interior of  $(\mu \times \Delta_{n-1})\cap \sigma$.
  The case when this hypothesis is violated can be modelled by
  adding additional linear constraints $f_i= 0$.
We can phrase our problem as a  parametric optimization problem:
    \begin{equation}\label{RosStuOP} 
        {\rm minimize}\,\,
    c_1x_1+\dots + c_nx_n \,\, \text{subject to}\,\,
\,    x \in \,\mathcal{M} =X\cap \Delta_{n-1}.
    \end{equation}
    Here $c_1,\ldots,c_n$ are parameters. In our applications, these $c_i$ 
    will be functions in the entries $d_{ij}$ of the metric $d$ and in the coordinates
 $\mu_k$ of the given point $\mu \in \Delta_{n-1}$. But, for now,  
 let us treat the $c_i$ as unknowns.
 The optimal value of the problem (\ref{RosStuOP}) is a function
 in these unknowns:
 $$ c_0^* \,\, = \,\, c_0^* (c_1,\ldots,c_n). $$
By \cite[Section 3]{RosStu10}, the {\em optimal value function}  
 $c_0^*:\mathbb{R}^n \rightarrow \mathbb{R}$ is an algebraic
 function in the $n$ parameters $c_1,\ldots ,c_n$. This means that
 there  exists a polynomial $\Phi(c_0,c_1,\dots, c_n)$ in $n+1$ variables such that 
 $\Phi(c_{0}^*,c_1,\dots ,c_n)$=0. The degree of $\Phi$ in its first argument $c_0$ 
 measures the algebraic complexity of our optimization problem~(\ref{RosStuOP}). 
 We call this number the {\em Wasserstein degree} of our model $\mathcal{M}$.
We shall describe the Wasserstein degree geometrically and offer some bounds.

Following \cite[Section 3]{RosStu10}, we consider the projective variety
$\bar{X}^*$ that is dual to the variety $\bar{X}$. The dual variety $\bar{X}^*$
lives in the dual projective space $\mathbb{P}^n$, and it parametrizes hyperplanes
in the ambient projective space of $\bar{X}$ that are tangent to $\bar{X}$.
This dual variety $\bar{X}^*$ is typically a hypersurface, regardless of
what the codimension of $X$ is. In particular, it is a hypersurface
when $X$ is compact in $\mathbb{R}^n$.  If $X$ is irreducible
then the hypersurface $\bar{X}^*$ is defined by a unique
(up to scaling) irreducible homogeneous polynomial in $n+1$ unknowns 
$c_0,c_1,\ldots,c_n$. The degree of this hypersurface is the
degree of $\bar{X}^*$. 
The following result is a direct consequence of  \cite[Theorem 3.2]{RosStu10}

\begin{theorem}\label{theoremRosStu} 
The polynomial $\Phi(-c_0,c_1,\ldots,c_n)$ is the defining equation
of the hypersurface $\bar{X}^*$ that is dual to the projective variety $\bar{X}$
that represents the  model $\mathcal{M}$ in $\Delta_{n-1}$.
Hence the Wasserstein degree of $\mathcal{M}$ is the degree
of $\Phi$ in its first argument. This is generically equal to the degree of $\bar{X}^*$.
\end{theorem}

For many natural classes of varieties $\bar{X}$, there are known
formulas for the degree of the dual $\bar{X}^*$.
This includes general complete intersections and determinantal varieties.
The case of a hypersurface appears in \cite[Example 2.7]{RosStu10}.
It serves as an illustration of our algebraic view on the problem (\ref{RosStuOP}).

\begin{corollary} \label{cor:hypersurface}
Suppose that the model $\mathcal{M}$ is a hypersurface,
namely, it is the zero
set in the simplex $\Delta_{n-1}$ of a general polynomial of degree $m$. Then
the Wasserstein degree of $\mathcal{M}$ equals~$m(m-1)^{n-2}$.
\end{corollary}

For instance, we have $n=m=2$ for the Hardy-Weinberg curve
(\ref{eq:hardyweinberg}), so this has Wasserstein degree $2$.
This reflects the fact that the optimal value 
$\sqrt{2} - 8/7$ is an algebraic number of degree~$2$.

\begin{example} If $\mathcal{M}$ is a general curve of degree $3$ in
the triangle $\Delta_2$ then its Wasserstein degree equals $6$.
Such an elliptic curve does not permit a rational parametrization,
so we will have to consider (\ref{RosStuOP}) as a constrained 
optimization problem.
For a concrete example consider the curve 
$$ x_1^3 + x_2^3 + x_3^3 \,\, = \,\, 4 x_1 x_2 x_3. $$
Let $c_0^*$ be the minimum of $c_1x_1+c_2x_2+c_3x_3$ over this curve in $\Delta_2$.
This is an algebraic function of degree $6$. Its minimal polynomial 
$\Phi(-c_0,c_1,c_2,c_3)$ is a homogeneous sextic. Namely, we have
$$ \begin{tiny} \begin{matrix}
\Phi \quad = \quad c_0^6 \,+\,
(2 c_1+2 c_2+2 c_3) c_0^5\,-\,
(65 c_1^2-70 c_1 c_2-70 c_1 c_3+65 c_2^2-70 c_2 c_3+65 c_3^2) c_0^4 \qquad
\qquad \qquad \\
+(208 c_1^3-442 c_1^2 c_2-442 c_1^2 c_3-442 c_1 c_2^2+2048 c_1 c_2 c_3-442 c_1 c_3^2+208 c_2^3-442 c_2^2 c_3-442 c_2 c_3^2+208 c_3^3) c_0^3 \\
-(117 c_1^4-546 c_1^3 c_2-546 c_1^3 c_3+1994 c_1^2 c_2^2-1024 c_1^2 c_2 c_3+1994 c_1^2 c_3^2 
-546 c_1 c_2^3-1024 c_1 c_2^2 c_3 \qquad \qquad \qquad \\ 
\qquad \qquad \qquad  -1024 c_1 c_2 c_3^2-546 c_1 c_3^3+117 c_2^4-546 c_2^3 c_3+
1994 c_2^2 c_3^2-546 c_2 c_3^3+117 c_3^4) c_0^2 \\
-(162 c_1^5-288 c_1^4 c_2-288 c_1^4 c_3+606 c_1^3 c_2^2-1152 c_1^3 c_2 c_3+606 c_1^3 c_3^2{+}
606 c_1^2 c_2^3{+}352 c_1^2 c_2^2 c_3{+}352 c_1^2 c_2 c_3^2{+}606 c_1^2 c_3^3
\\ -288 c_1 c_2^4 -1152
c_1 c_2^3 c_3+352 c_1 c_2^2 c_3^2-1152 c_1 c_2 c_3^3-288 c_1 c_3^4{+}162 c_2^5-288 c_2^4 c_3{+}
606 c_2^3 c_3^2{+}606 c_2^2 c_3^3 {-}288 c_2 c_3^4 \\ +162 c_3^5) c_0 
-27 c_1^6+288 c_1^4 c_2 c_3{-}202 c_1^3 c_2^3{-}202 c_1^3 c_3^3{-}176 c_1^2 c_2^2 c_3^2
{+}288 c_1c_2^4 c_3{+}288 c_1 c_2 c_3^4{-}27 c_2^6{-}202 c_2^3 c_3^3{-}27 c_3^6.
\end{matrix}
\end{tiny}
$$
For any given $c_1,c_2,c_3$, the optimal value is obtained by
solving $\Phi = 0$ for $c_0$.
\end{example}

As we said earlier, in our application in Step 5 of Algorithm \ref{algorithm},
the $c_i$ depend on the matrix $d$ and the distribution $\mu$.
We can thus consider the function that measures the Wasserstein distance:
    $$ \mathbb{R}^{n^2}\times \Delta_{n-1} \,\rightarrow \,\mathbb{R}\,,\quad
    (d,\mu) \,\mapsto \,c_{0}^*(d,\mu)\,=\,W_d(\mu,\mathcal{M}). $$
Our discussion establishes the following result about this function
which depends only on $\mathcal{M}$.

\begin{corollary} \label{cor:piece}
The Wasserstein distance is a piecewise algebraic function of $d$ and $\mu$.
Each piece is an algebraic function whose degree is bounded above
by the degree of the hypersurface dual to $\mathcal{M}$.
\end{corollary}

\section{The Wasserstein Estimator of an Independence Model}
\label{sec:independence}

In this section we present our solution to the problem that started this project.
The task is to compute the Wasserstein estimator for the independence model
on two binary random variables. Here $n=4$ and $\mathcal{M}$
is the variety of $2 \times 2$ matrices of rank $1$. This has the
parametric representation
$$ \quad \begin{pmatrix} x_1 & x_2 \\ x_3 & x_4  \end{pmatrix}
\,\, = \,\,
\begin{pmatrix} pq & \,\,p(1-q) \\ (1{-}p)q & \,\,(1{-}p)(1{-}q) \end{pmatrix},
\,\, \text{where}\,\,  (p,q)\in[0,1]^2. 
$$
Equivalently, $\mathcal{M}$ is the quadratic surface $  \{x_1 x_4 = x_2 x_3 \}$
in the tetrahedron $\Delta_3$.

Our underlying metric space $\Omega $ is the square $\{0,1\}^2$ with its Hamming distance.
We identify $\Omega$ with the set $[4] = \{1,2,3,4\}$
as indicated above. The ground metric is 
represented by the matrix
$$
d \,\,=\,\,\,
\begin{blockarray}{*{4}{c} l}
    \begin{block}{*{4}{>{\small}c<{}} l}
     1 & 2 & 3 & 4  & \\
    \end{block}
     \begin{block}{[*{4}{c}]>{\;\;\small}l<{}}
\,0 & 1 & 1 & 2 \,\, &  1  \\
\,1 & 0 & 2 & 1  \,\, & 2  \\
\,1 & 2 & 0 & 1 \,\,& 3  \\
\,2 & 1 & 1 & 0\,\, & 4 \\
\end{block} 
\end{blockarray}\vspace{-0.15in}
$$ 
Given two points $\mu,\nu$ in  $\Delta_3$, the
transportation polytope $\Pi(\mu,\nu)$
 consists of all nonnegative $4 \times 4$ matrices $\pi$ with row sums  
$\nu$ and column sums $\mu$. 
It usually is simple and has dimension $9$.
The Wasserstein distance between the two distributions equals
$W(\mu,\nu) = \min_{\pi\in \Pi(\mu,\nu)} \sum_{1\leq i,j\leq 4} d_{ij}\pi_{ij}$. 

What we are interested in is the minimum Wasserstein distance from $\mu$
to any point $\nu$ in the independence model $\mathcal{M}$. 
This parametric linear optimization problem can be described as follows:
\[
\begin{blockarray}{ccccc}
\begin{block}{[cccc]l}
  \pi_{11} & \pi_{12} & \pi_{13} & \pi_{14} & \;\; pq \\
  \pi_{21} & \pi_{22} & \pi_{23} & \pi_{24} & \;\; p(1-q) \\
  \pi_{31} & \pi_{32} & \pi_{33} & \pi_{34} & \;\; (1-p)q \\
  \pi_{41} & \pi_{42} & \pi_{43} & \pi_{44} & \;\; (1-p)(1-q) \\
\end{block}
\begin{block}{ccccc}
\mu_1 & \mu_2 & \mu_3 & \mu_4 & \\
\end{block}
\end{blockarray}   \vspace{-0.15in}
 \]
Here the marginal $\mu=(\mu_1,\mu_2,\mu_3,\mu_4)$ is fixed.
The model $\mathcal{M}$ is parametrized by
the points $(p,q)$ in the square $[0,1]^2$. 
The Wasserstein distance between $\mu$ and  $\nu = \nu(p,q)$
is a continuous function on that square. The minimum value of that
function is the desired Wasserstein distance
$W(\mu,\mathcal{M})$.

\begin{figure}[h]
    \centering
    \vspace{-0.15in}
\includegraphics[scale=.45,clip=true, trim=0cm 0.5cm 15cm 3cm]{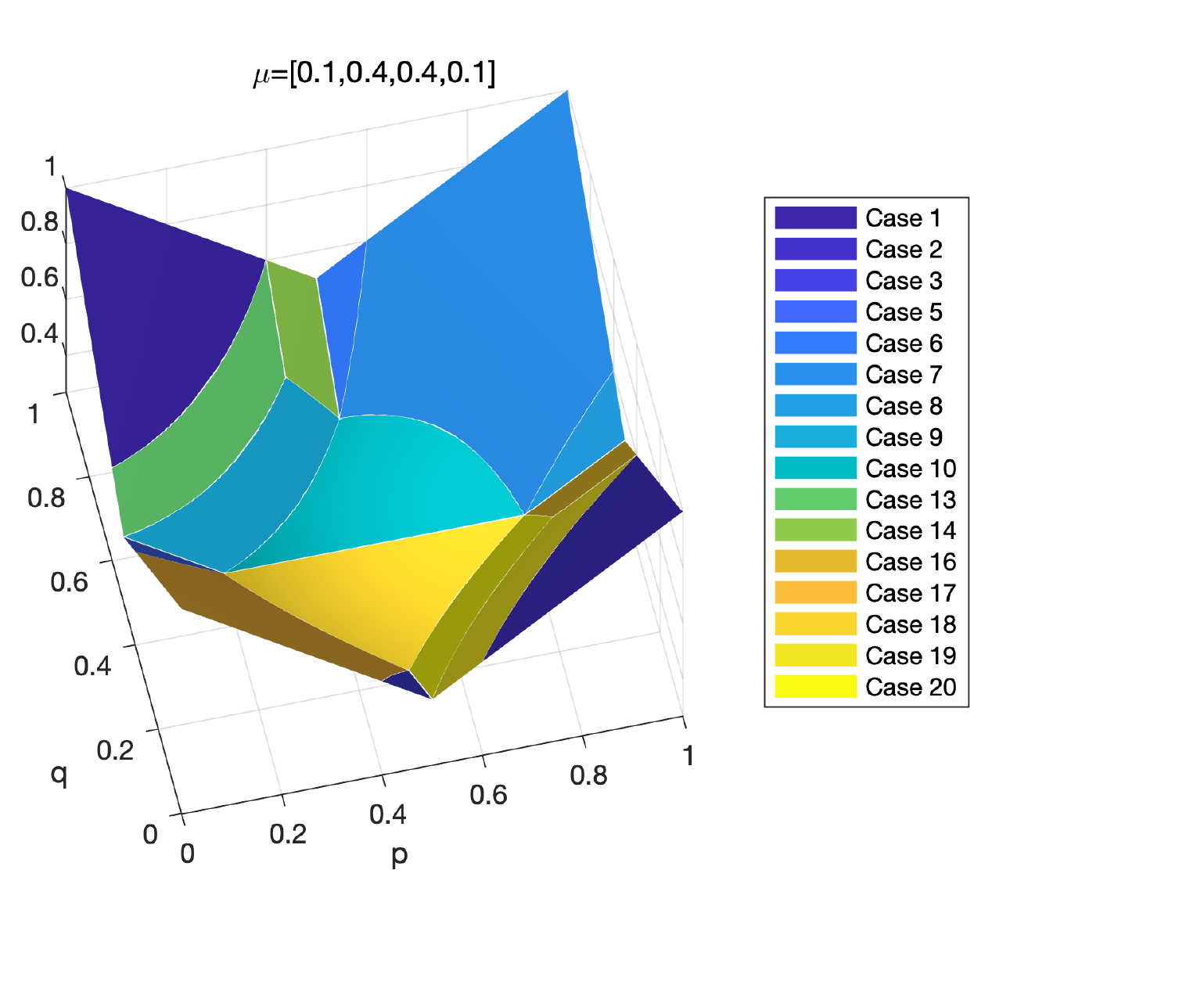} \!\! \!
\includegraphics[scale=.32,clip=true, trim=15cm -3cm 0cm 0cm]{target01_04_04_01_1}
\vspace{-0.38in}
    \caption{Left: The graph of the distance function $\,[0,1]^2 \rightarrow \mathbb{R},\,
    (p,q) \mapsto W(\mu,\nu(p,q))\,$
    for $\,\mu = \frac{1}{10} (1,4,4,1)$. Right: The independence model $\mathcal{M}$ inside the 
    tetrahedron $\Delta_3$. Color corresponds to the Wasserstein distance to the target distribution $\mu$, shown as a red dot. The function has two global minimizers over $\mathcal{M}$.
    \label{fig:indep_mod1}}
\end{figure}    

Our task is to minimize the function in Figure \ref{fig:indep_mod1} over
the square. We see that this function is piecewise algebraic 
(cf.~Corollary \ref{cor:piece}). Each piece is either a linear
function or a quadratic function.  The case distinction arises from
the induced polyhedral subdivision of the $6$-dimensional
polytope $\Delta_3 \times \Delta_3$. This subdivision is not a triangulation,
but, following Step 2 in Algorithm \ref{algorithm}, we replace
it with a nearby triangulation. That triangulation has $20$
maximal simplices, as seen in (\ref{eq:stirling}). These
are the $20$ cases in Figure \ref{fig:indep_mod1}.
The graph of our function is color-coded according to these cases.

\begin{figure}[t]
\centering
\includegraphics[scale=.35,clip=true]{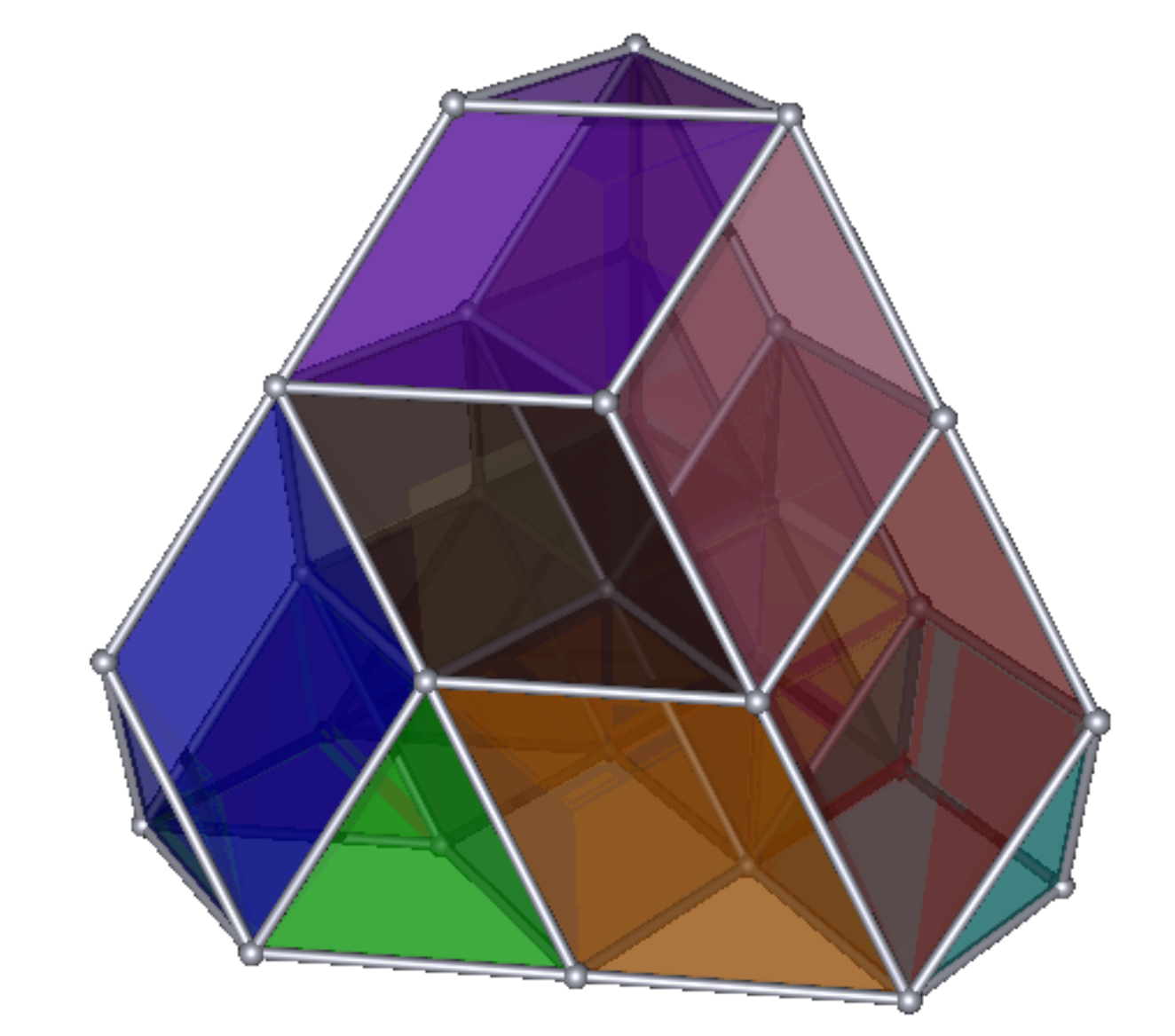} \qquad \quad
\includegraphics[scale=.3,clip=true]{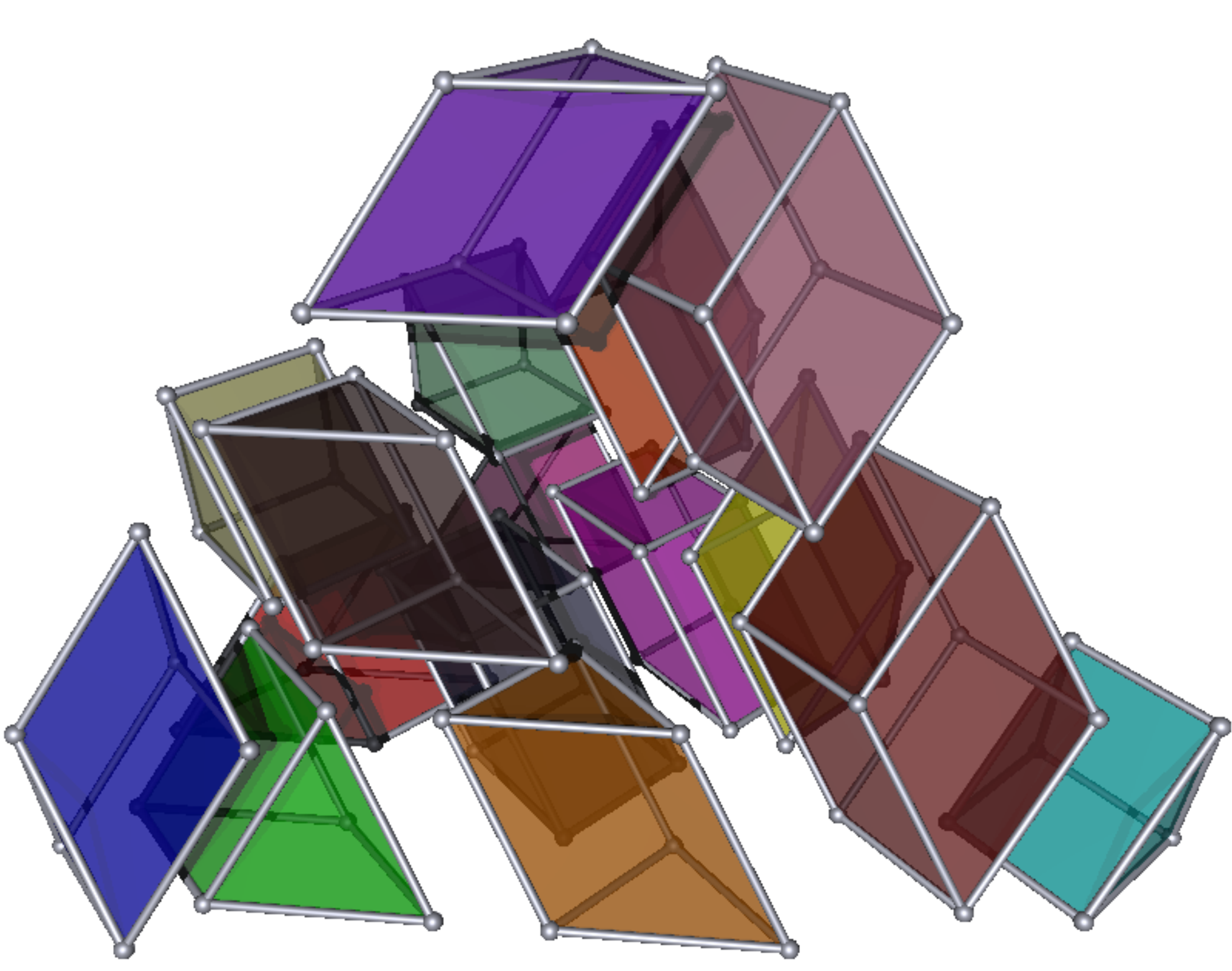}
\caption{A mixed subdivision of a truncated tetrahedron
into $16=12+4$ cells. \label{fig:explosion}}
\end{figure}

The triangulation of $\Delta_3 \times \Delta_3$
restricts to a mixed subdivision of the tetrahedron
$\mu \times \Delta_3$. That subdivision consists of 
$20 = 4+12+4$ cells, namely $4$ tetrahedra,
$12$ triangular prisms, and $4$ parallelepipeds. 
After removing the $4$ tetrahedra, which 
touch the vertices of $\mu \times \Delta_3$, we obtain a
{\em truncated tetrahedron} which is subdivided into $16$ cells.
Such a subdivision is shown in Figure~\ref{fig:explosion}.

The restriction of the mixed subdivision of $\mu \times \Delta_3$ divides our model $\mathcal{M}$ into regions. On each of these regions, 
the Wasserstein distance function $\nu \mapsto W(\mu,\nu)$
is given by a linear functional, as explained 
in Step 4 of Algorithm \ref{algorithm}. The surface
and this function are depicted in Figure \ref{fig:indep_mod1} (right).


The two images in \Cref{fig:indep_mod1} convey the same information.
The piecewise linear function on the quadratic surface in $\Delta_3$
restricts to a piecewise quadratic function on the square $[0,1]^2$
under the parametrization of the surface. However, the color
coding in the two diagrams is different. The colors in the left image in Figure~\ref{fig:indep_mod1} 
show the pieces, while the right one displays a heat map.
Namely, the colors here  represent values of the function
$\mathcal{M} \rightarrow \mathbb{R},\,\nu \mapsto W(\mu,\nu)$.
The two bluest points attain the minimum value $W(\mu,\mathcal{M})$.
The white curve segments on the surface $\mathcal{M}$
are the boundaries between the various pieces. 
Each piece is the intersection of $\mathcal{M}$ with one of the
polytopes in the mixed subdivision  in Figure~\ref{fig:explosion}.

We now discuss the computations that led to these results and pictures.
The triangulation of $\Delta_3 \times \Delta_3$ and
the resulting mixed subdivision  of $\mu \times \Delta_3$
are computed in Steps 1-3 of Algorithm~\ref{algorithm}.
These geometric objects are encoded algebraically, namely in the
decomposition (\ref{eq:primdec}) of the ideal
$$
\begin{scriptsize}
\begin{matrix}
\text{in}_{d_{\epsilon}}(I_A)\,\,=\,\,&\langle y_{11},y_{12},y_{14},y_{31},y_{32},y_{34},y_{41},y_{42},y_{44}\rangle \,\,\,\cap &
&\langle y_{11},y_{13},y_{14},y_{21},y_{23},y_{24},y_{41},y_{43},y_{44}\rangle \,\,\,\cap\\
&\langle y_{12},y_{1  
	3},y_{14},y_{23},y_{24},y_{31},y_{32},y_{34},y_{42}\rangle \,\,\,\cap &
&\langle y_{12},y_{13},y_{14},y_{2    
	3},y_{24},y_{32},y_{34},y_{42},y_{43}\rangle \,\,\,\cap \\
& \langle y_{12},y_{13},y_{14},y_{21},y_{23},y_{24},y_{32},y_{34},y_{43}\rangle \,\,\,\cap   &
&\langle y_{13},y_{21},y_{23},y_{24},y_{31},y_{32},y_{41},y_{42},y_{43}\rangle \,\,\,\cap          \\
& \langle y_{21},y_{2  
	3},y_{24},y_{31},y_{32},y_{34},y_{41},y_{42},y_{43}\rangle\,\,\,\cap &
&\langle y_{12},y_{21},y_{23},y_{3     
	1},y_{32},y_{34},y_{41},y_{42},y_{43}\rangle \,\,\,\cap \\
&\langle y_{13},y_{14},y_{21},y_{23},y_{24},y_{3    
	2},y_{34},y_{41},y_{43}\rangle \,\,\,\cap  &
&\langle y_{14},y_{21},y_{23},y_{24},y_{31},y_{32},y_{34},y_{41},y_{43}\rangle \,\,\,\cap\\
&\langle y_{12},y_{13},y_{21},y_{23},y_{31},y_{32},y_{41},y_{4     
	2},y_{43}\rangle \,\,\,\cap &
&\langle y_{12},y_{13},y_{14},y_{21},y_{23},y_{32},y_{41},y_{42},y_{43}\rangle \,\,\,\cap\\
& \langle y_{11},y_{13},y_{14},y_{2     
	1},y_{23},y_{24},y_{32},y_{41},y_{43}\rangle \,\,\,\cap &
&\langle y_{13},y_{14},y_{21},y_{23},y_{24},y_{32},y_{41},y_{4     
	2},y_{43}\rangle \,\,\,\cap \\
&\langle y_{12},y_{13},y_{14},y_{23},y_{32},y_{3    
	4},y_{41},y_{42},y_{43}\rangle \,\,\,\cap &
&\langle y_{12},y_{1  
	4},y_{23},y_{31},y_{32},y_{34},y_{41},y_{42},y_{43}\rangle \,\,\,\cap \\
& \langle y_{12},y_{13},y_{14},y_{21},y_{23},y_{24},y_{31},y_{3    
	2},y_{34}\rangle \,\,\,\cap &
&\langle y_{12},y_{14},y_{21},y_{23},y_{24},y_{31},y_{32},y_{34},y_{41}\rangle \,\,\,\cap     \\     
&\langle y_{11},y_{12},y_{14},y_{23},y_{31},y_{32},y_{34},y_{41},y_{42}\rangle \,\,\,\cap &
& \!\!\!\!\! \langle y_{12},y_{14},y_{23},y_{24},y_{31},y_{32},y_{34},y_{41},y_{42}\rangle .
\end{matrix} \end{scriptsize}
$$
Step 4 of Algorithm~\ref{algorithm} translates each of these $20$ minimal primes
into a $4 \times 4$-matrix in the variety of that prime 
whose nonzero entries are linear forms in $\mu_i$ and $\nu_j$. For instance, the
second prime~gives
\begin{equation}
\label{eq:census} \tilde \pi_\sigma \,\, = \,\, 
\begin{small} \begin{bmatrix}
\,\, 0 & \nu_1 & 0 & \,0 \,\,\\
\,\, 0 & \nu_2 & 0 & \,0 \,\,\\
\,\,\mu_1\,\, &\, \mu_2 {-} \nu_1{-}\nu_2{-}\nu_4 \,&\,\,\, \mu_3 \,\,\,& \,\mu_4\,\, \\ 
\,\, 0 & \nu_4 & 0 & \,0 \,\,
\end{bmatrix} \end{small}.
\end{equation}
The dot product of  $d$ and $\tilde \pi_\sigma$ gives
the Wasserstein distance on the piece labeled {\em Case 2} in Figure~\ref{fig:indep_mod1}:
$$ d \cdot \tilde \pi_\sigma \,\, = 
\,\,
\mu_1 {+} 2 \mu_2 {+} \mu_4 - \nu_1 {-} 2 \nu_2 {-} \nu_4 \,\, = \,\,
        \mu_2 - \mu_3 -\nu_2 + \nu_3\,\, = \,\, \mu_2 - \mu_3 - p + q.
$$
Hence {\em Case 2} is linear in $p,q$. 
The region in the square for this case is defined by the requirement 
that the entries of the matrix $\tilde \pi_\sigma$ are between $0$ and $1$. 
We only need to consider the entry in the third row and second column:
$$ 0 \,\,\, \leq \,\,\, \mu_2 -  \nu_1 - \nu_2 - \nu_4 \,\, =\,\,
    (1-p)q+\mu_2-1 \,\,\, \leq \,\,\, 1 . $$
In Figure~\ref{fig:indep_mod1}
the graph of $\mu_2 - \mu_3 - p + q$ on this region is shown in blue and labeled {\em Case 2}. 

We analyze all $20$ components of ${\rm in}_d(I_A)$ in this manner,
and we record the result in the first two columns of Tables \ref{tab:bigg1} and \ref{tab:bigg2}. 
The rightmost column gives the support of the corresponding
vertex of the transportation polytope. For instance, the matrix in the first row of the table shows the support of $\tilde \pi_\sigma$ in (\ref{eq:census}). 

The third column of the table contains all candidates for the optimal point $\nu^*$ expressed as an algebraic function in the four coordinates of $\mu$. 
Each of the $20$ cases has one or several candidate solutions, listed in the third and forth columns of the table. 
Which of the candidates is the actual solution can be determined in terms of further case distinctions on $\mu$, which we omit in the table. 
%
%
The smallest solution among all cases for a given $\mu$ is the desired Wasserstein distance $W(\mu,\nu^*) = W(\mu,\mathcal{M})$. 
Note that these expressions involve a square root, so the Wasserstein degree of the independence surface $\mathcal{M}$ equals two, as predicted by Corollary \ref{cor:hypersurface}. 

\subsection*{Conclusion} In this paper, we developed mathematical foundations for computing the Wasserstein distance between a point and an algebraic variety in a probability simplex. Our next goal is to develop a practical algorithm that scales beyond toy problems. We also plan to answer the questions raised in the introduction, such as characterizing scenarios when the minimizer is unique.

\newcommand{\tone}{%
\left[\begin{smallmatrix}
0 & * & 0 & 0 \\
0 & * & 0 & 0 \\
* & * & * & * \\
0 & * & 0 & 0 \\
\end{smallmatrix}\right]
}
\newcommand{\ttwo}{%
\left[\begin{smallmatrix}
0 & 0 & * & 0 \\
* & * & * & * \\
0 & 0 & * & 0 \\
0 & 0 & * & 0 \\
\end{smallmatrix}\right]
}
\newcommand{\tthree}{%
\left[\begin{smallmatrix}
* & 0 & 0 & 0 \\
* & * & 0 & 0 \\
0 & 0 & * & 0 \\
* & 0 & * & * \\
\end{smallmatrix}\right]
}
\newcommand{\tfour}{%
\left[\begin{smallmatrix}
* & 0 & 0 & 0 \\
* & * & 0 & 0 \\
* & 0 & * & 0 \\
* & 0 & 0 & * \\
\end{smallmatrix}\right]
}
\newcommand{\tfive}{%
\left[\begin{smallmatrix}
* & 0 & 0 & 0 \\
0 & * & 0 & 0 \\
* & 0 & * & 0 \\
* & * & 0 & * \\
\end{smallmatrix}\right]
}
\newcommand{\tsix}{%
\left[\begin{smallmatrix}
* & * & 0 & * \\
0 & * & 0 & 0 \\
0 & 0 & * & * \\
0 & 0 & 0 & * \\
\end{smallmatrix}\right]
}
\newcommand{\tseven}{%
\left[\begin{smallmatrix}
* & * & * & * \\
0 & * & 0 & 0 \\
0 & 0 & * & 0 \\
0 & 0 & 0 & * \\
\end{smallmatrix}\right]
}
\newcommand{\teight}{%
\left[\begin{smallmatrix}
* & 0 & * & * \\
0 & * & 0 & * \\
0 & 0 & * & 0 \\
0 & 0 & 0 & * \\
\end{smallmatrix}\right]
}
\newcommand{\tnine}{%
\left[\begin{smallmatrix}
* & * & 0 & 0 \\
0 & * & 0 & 0 \\
* & 0 & * & 0 \\
0 & * & 0 & * \\
\end{smallmatrix}\right]
}
\newcommand{\tten}{%
\left[\begin{smallmatrix}
* & * & * & 0 \\
0 & * & 0 & 0 \\
0 & 0 & * & 0 \\
0 & * & 0 & * \\
\end{smallmatrix}\right]
}
\newcommand{\televen}{%
\left[\begin{smallmatrix}
* & 0 & 0 & * \\
0 & * & 0 & * \\
0 & 0 & * & * \\
0 & 0 & 0 & * \\
\end{smallmatrix}\right]
}
\newcommand{\ttwelve}{%
\left[\begin{smallmatrix}
* & 0 & 0 & 0 \\
0 & * & 0 & * \\
* & 0 & * & * \\
0 & 0 & 0 & * \\
\end{smallmatrix}\right]
}
\newcommand{\tthirteen}{%
\left[\begin{smallmatrix}
0 & * & 0 & 0 \\
0 & * & 0 & 0 \\
* & 0 & * & * \\
0 & * & 0 & * \\
\end{smallmatrix}\right]
}
\newcommand{\tfourteen}{%
\left[\begin{smallmatrix}
* & * & 0 & 0 \\
0 & * & 0 & 0 \\
* & 0 & * & * \\
0 & 0 & 0 & * \\
\end{smallmatrix}\right]
}
\newcommand{\tfifteen}{%
\left[\begin{smallmatrix}
* & 0 & 0 & 0 \\
* & * & 0 & * \\
* & 0 & * & 0 \\
0 & 0 & 0 & * \\
\end{smallmatrix}\right]
}
\newcommand{\tsixteen}{%
\left[\begin{smallmatrix}
* & 0 & * & 0 \\
* & * & 0 & * \\
0 & 0 & * & 0 \\
0 & 0 & 0 & * \\
\end{smallmatrix}\right]
}
\newcommand{\tseventeen}{%
\left[\begin{smallmatrix}
* & 0 & 0 & 0 \\
0 & * & 0 & 0 \\
0 & 0 & * & 0 \\
* & * & * & * \\
\end{smallmatrix}\right]
}
\newcommand{\teighteen}{%
\left[\begin{smallmatrix}
* & 0 & * & 0 \\
0 & * & 0 & 0 \\
0 & 0 & * & 0 \\
0 & * & * & * \\
\end{smallmatrix}\right]
}
\newcommand{\tnineteen}{%
\left[\begin{smallmatrix}
0 & 0 & * & 0 \\
* & * & 0 & * \\
0 & 0 & * & 0 \\
0 & 0 & * & * \\
\end{smallmatrix}\right]
}
\newcommand{\ttwenty}{%
\left[\begin{smallmatrix}
* & 0 & * & 0 \\
* & * & 0 & 0 \\
0 & 0 & * & 0 \\
0 & 0 & * & * \\
\end{smallmatrix}\right]
}

\newcolumntype{L}{>{$}c<{$}}  

\begin{table}[h!]
\begin{center}
\vspace{-0.14in}
\scalebox{.7}{
  \begin{tabular}{|c | L | L | L | L | L | }
  \hline
Case &   \text{Objective Function} & \text{Feasible Region, $0\leq *\leq 1$} & \text{Solution} & \text{Minimum Value} & \text{Subdivision} \\\hline\hline                     
\multicolumn{6}{|l|}{\quad Quadratic pieces}\\\hline 
10& 2 p q-p-q+\mu_2+\mu_3&\begin{tabular}{L}
                         q-\mu_1-\mu_3\\
                         \mu_3-(1-p) q\\
                         (1-p) (1-q)-\mu_4
                        \end{tabular}&
                        \begin{tabular}{L}
                         (\frac{\mu_1}{\mu_1+\mu_3},\mu_1+\mu_3)\\
                         (\frac{\mu_2}{\mu_2+\mu_4},\mu_1+\mu_3)\\
                         (\mu_1+\mu_2,\frac{\mu_3}{\mu_3+\mu_4})\\
                         \left(1-\sqrt{\mu_3},\sqrt{\mu_3}\right)
                        \end{tabular}
                        &
                        \begin{tabular}{L}
                        -\frac{\mu_1}{\mu_1+\mu_3}+\mu_1+\mu_2 \\
                        \frac{\mu_2}{\mu_2+\mu_4}-\mu_1-\mu_2\\
                        \frac{\mu_3}{\mu_3+\mu_4}-\mu_1-\mu_3\\ 
                        2\sqrt{\mu_3}(1-\sqrt{\mu_3})-\mu_1-\mu_4
                        \end{tabular}& \tten \\\hline
18& 2 p q-p-q+\mu_2+\mu_3&\begin{tabular}{L}
                        p q-\mu_1\\
                        \mu_2- p (1-q)\\
                        \mu_1+\mu_3-q
                     \end{tabular}&
                     \begin{tabular}{L}
                       \left(\mu_1+\mu_2,\frac{\mu_1}{\mu_1+\mu_2}\right)\\
                       \left(\frac{\mu_2}{\mu_2+\mu_4},\mu_1+\mu_3\right)\\
                       (\frac{\mu_1}{\mu_1+\mu_3},\mu_1+\mu_3)\\
                       (\sqrt{\mu_2},1-\sqrt{\mu_2})
                     \end{tabular}&
                     \begin{tabular}{L}
                       -\frac{\mu_1}{\mu_1+\mu_2}+\mu_1+\mu_3\\
                      \frac{\mu_2}{\mu_2+\mu_4}- \mu_1-\mu_2\\
                       -\frac{\mu_1}{\mu_1+\mu_3}+\mu_1+\mu_2\\
                       2\sqrt{\mu_2}(1-\sqrt{\mu_2})-\mu_1-\mu_4
                     \end{tabular}&  \teighteen \\\hline
12& -2 p q +p+q-\mu_2-\mu_3&\begin{tabular}{L}
                    \mu_1-p q\\
                    p (1-q)-\mu_2\\
                     q -\mu_1-\mu_3
                     \end{tabular}&
                     \begin{tabular}{L}
                     (\mu_1+\mu_2,\frac{\mu_1}{\mu_1+\mu_2})\\
                     (\frac{\mu_1}{\mu_1+\mu_3},\mu_1+\mu_3)\\
                     (\frac{\mu_2}{\mu_2+\mu_4},\mu_1+\mu_3)\\
                     (\sqrt{\mu_1},\sqrt{\mu_1})
                     \end{tabular}
                     & 
                     \begin{tabular}{L}
                       \frac{\mu_1}{\mu_1+\mu_2}-\mu_1-\mu_3\\                       
                        \frac{\mu_1}{\mu_1+\mu_3}-\mu_1-\mu_2\\
                        -\frac{\mu_2}{\mu_2+\mu_4}+\mu_1+\mu_2 \\
                        2\sqrt{\mu_1}(1-\sqrt{\mu_1})-\mu_2-\mu_3
                     \end{tabular}& \ttwelve \\\hline
15& -2 p q+p+q-\mu_2-\mu_3&\begin{tabular}{L}
                    \mu_1+\mu_3-q\\
                    (1-p) q -\mu_3\\
                    \mu_4-(1-p) (1-q)
                     \end{tabular}& 
                     \begin{tabular}{L}
                       \left(1-\sqrt{\mu_4},1-\sqrt{\mu_4}\right)\\
                       \left(\frac{\mu_1}{\mu_1+\mu_3}, \mu_1+\mu_3\right)\\
                      \left(\mu_1+\mu_2,\frac{\mu_3}{\mu_3+\mu_4}\right)\\
                      \left(\frac{\mu_2}{\mu_2+\mu_4},\mu_1+\mu_3 \right)
                     \end{tabular}&
                     \begin{tabular}{L}
                     2\sqrt{\mu_4}(1-\sqrt{\mu_4})-\mu_2-\mu_3\\
                       \frac{\mu_1}{\mu_1+\mu_3}-\mu_1-\mu_2\\
                      -\frac{\mu_3}{\mu_3+\mu_4}+\mu_1+\mu_3\\
                      -\frac{\mu_2}{\mu_2+\mu_4} +\mu_1+\mu_2                      
                     \end{tabular} & \tfifteen \\\hline
  \end{tabular}
}  
\medskip
\caption{Algebraic analysis of the Wasserstein distance function shown
in Figure \ref{fig:indep_mod1}.
 \label{tab:bigg2}}
  \vspace{-0.2in}
  \end{center}
  \end{table}

\begin{table}[h!]
\begin{center}
\vspace{-0.14in}
\scalebox{.7}{
  \begin{tabular}{|c | L | L | L | L | L | }
  \hline
Case &   \text{Objective Function} & \text{Feasible Region, $0\leq *\leq 1$} & \text{Solution} & \text{Minimum Value} & \text{Subdivision} \\\hline\hline  
\multicolumn{6}{|l|}{\quad First affine piece}\\\hline
2&  -p+q+\mu_2-\mu_3  & (1-p)q+\mu_2-1 & (1-\sqrt{1-\mu_2},\sqrt{1-\mu_2})
 & 2\sqrt{1-\mu_2}+\mu_2-\mu_3-1 & \tone \\[10pt] \hline
9& -p+q+\mu_2-\mu_3&\begin{tabular}{L}
                    \mu_1+\mu_3-(1-p) q \\
                    q-\mu_1-\mu_3\\
                    (1-p) q -\mu_3\\
                    (1-p) (1-q)-\mu_4
                     \end{tabular}& \begin{tabular}{L}
                      \left(\frac{\mu_2}{\mu_2+\mu_4},\mu_1+\mu_3\right) \\
                     (1-\sqrt{\mu_3},\sqrt{\mu_3}) \\
                     \left(\frac{\mu_1}{\mu_1+\mu_3},\mu_1+\mu_3\right) \\
                      \left(\mu_1+\mu_2,\frac{\mu_3}{\mu_3+\mu_4}\right)
                     \end{tabular} & 
                      \begin{tabular}{L}
                      -\frac{\mu_2}{\mu_2+\mu_4} +\mu_1+\mu_2  \\
                     2\sqrt{\mu_3}+ \mu_2-\mu_3-1 \\
                    -\frac{\mu_1}{\mu_1+\mu_3} +\mu_1+\mu_2 \\
                     \frac{\mu_3}{\mu_3+\mu_4} - \mu_1-\mu_3
                     \end{tabular} &  \tnine \\\hline
13& -p+q+\mu_2-\mu_3&\begin{tabular}{L}
                    \mu_2-p\\
                    (1-p) q -\mu_1-\mu_3\\
                    1-\mu_2-(1-p) q
                     \end{tabular}& \begin{tabular}{L}
                      \left(1-\sqrt{\mu_1+\mu_3},\sqrt{\mu_1+\mu_3}\right)\\
                     \left(\mu_2,\frac{\mu_1+\mu_3}{\mu_1+\mu_3+\mu_4}\right)
                     \end{tabular} & 
                      \begin{tabular}{L}
                      2\sqrt{\mu_1+\mu_3}+\mu_2-\mu_3-1  \\
                     \frac{\mu_1+\mu_3}{\mu_1+\mu_3+\mu_4} -\mu_3 
                     \end{tabular} & \tthirteen \\\hline
14& -p+q+\mu_2-\mu_3&\begin{tabular}{L}
                    p-\mu_2\\
                    \mu_2-p (1-q)\\
                    \mu_1+\mu_2-p\\
                    \mu_4-(1-p) (1-q)
                     \end{tabular}& \begin{tabular}{L}
                      \left(\mu_1+\mu_2,\frac{\mu_3}{\mu_3+\mu_4}\right)\\
                     \left(\mu_1+\mu_2,\frac{\mu_1}{\mu_1+\mu_2}\right)\\
                     \left(\frac{\mu_2}{\mu_2+\mu_4},\mu_1+\mu_3\right)
                     \end{tabular} & 
                      \begin{tabular}{L}
                      \frac{\mu_3}{\mu_3+\mu_4}-\mu_1-\mu_3 \\
                     \frac{\mu_1}{\mu_1+\mu_2}-\mu_1-\mu_3 \\
                     -\frac{\mu_2}{\mu_2+\mu_4}+\mu_1+\mu_2
                     \end{tabular} & \tfourteen \\\hline\hline 

\multicolumn{6}{|l|}{\quad Second affine piece}\\\hline                     
1&  p-q-\mu_2+\mu_3 &  p (1-q)+\mu_3-1 & (\sqrt{1-\mu_3},1-\sqrt{1-\mu_3}) &  2\sqrt{1-\mu_3}-\mu_2+\mu_3-1 & \ttwo \\[10pt]\hline
16& p-q-\mu_2+\mu_3&\begin{tabular}{L}
                    q-\mu_3\\
                    \mu_3-(1-p) q \\
                    \mu_1+\mu_3-q\\
                    \mu_4-(1-p) (1-q)
                     \end{tabular}& \begin{tabular}{L} \left(\mu_1+\mu_2,\frac{\mu_3}{\mu_3+\mu_4}\right)\\\left(\frac{\mu_1}{\mu_1+\mu_3},\mu_1+\mu_3\right)\\   \left(\frac{\mu_2}{\mu_2+\mu_4},\mu_1+\mu_3\right)
                     \end{tabular}& 
                   \begin{tabular}{L} -\frac{\mu_3}{\mu_3+\mu_4}+\mu_1+\mu_3 \\ \frac{\mu_1}{\mu_1+\mu_3}-\mu_1-\mu_2\\
                   \frac{\mu_2}{\mu_2+\mu_4} -\mu_1-\mu_2 
                     \end{tabular} & \tsixteen \\\hline
19& p-q-\mu_2+\mu_3&\begin{tabular}{L}
                    p (1-q)-\mu_1-\mu_2\\
                    \mu_3- q \\
                    1-\mu_3-p (1-q)
                     \end{tabular}&
                     \begin{tabular}{L}
                   \left(\sqrt{\mu_1+\mu_2},1-\sqrt{\mu_1+\mu_2}\right)\\
                   \left(\frac{\mu_1+\mu_2}{\mu_1+\mu_2+\mu_4},\mu_3\right)
                     \end{tabular}&
                    \begin{tabular}{L}
                   2\sqrt{\mu_1+\mu_2}-\mu_2+\mu_3-1\\
                   \frac{\mu_1+\mu_2}{\mu_1+\mu_2+\mu_4}-\mu_2
                     \end{tabular}& \tnineteen \\\hline
20& p-q-\mu_2+\mu_3&\begin{tabular}{L}
                    \mu_1+\mu_2-p (1-q)\\
                     p (1-q)-\mu_2\\
                    p-\mu_1-\mu_2\\
                    (1-p) (1-q)-\mu_4
                     \end{tabular}&
                     \begin{tabular}{L}
                   \left(\mu_1+\mu_2,\frac{\mu_3}{\mu_3+\mu_4}\right)\\
                   \left(\sqrt{\mu_2},1-\sqrt{\mu_2}\right)\\
                   \left(\mu_1+\mu_2,\frac{\mu_1}{\mu_1+\mu_2},\right)\\
                   \left(\frac{\mu_2}{\mu_2+\mu_4},\mu_1+\mu_3\right)
                     \end{tabular}&
                    \begin{tabular}{L}
                    -\frac{\mu_3}{\mu_3+\mu_4}+\mu_1+\mu_3\\
                    2\sqrt{\mu_2}-\mu_2+\mu_3-1\\
                    -\frac{\mu_1}{\mu_1+\mu_2}+\mu_1+\mu_3\\
                    \frac{\mu_2}{\mu_2+\mu_4}-\mu_1-\mu_2
                     \end{tabular} 
                     &  \ttwenty \\ \hline\hline
\multicolumn{6}{|l|}{\quad Third affine piece}\\\hline 
3& -p-q+\mu_1-\mu_4+1 & \begin{tabular}{L} 
                       p (1-q)-\mu_2 \\
                       \mu_1+\mu_2-p \\
                       \mu_3-(1-p) q \end{tabular}  &
                       \begin{tabular}{L} 
                       \left(\mu_1+\mu_2,\frac{\mu_1}{\mu_1+\mu_2}\right)\\
                       \left(\mu_1+\mu_2,\frac{\mu_3}{\mu_3+\mu_4}\right)
                       \end{tabular}
                       &
                       \begin{tabular}{L} 
                       -\frac{\mu_1}{\mu_1+\mu_2}+\mu_1+\mu_3\\
                       -\frac{\mu_3}{\mu_3+\mu_4}+\mu_1+\mu_3
                       \end{tabular}
                       & \tthree \\\hline
4& -p-q+\mu_1-\mu_4+1 & \begin{tabular}{L}
                      p (1-q)-\mu_2\\
                      (1-p) q-\mu_3\\
                      (1-p) (1-q)-\mu_4
                      \end{tabular}&
                      \begin{tabular}{L}
                      \left(\gamma^+,1-\frac{\mu_2}{\gamma^+}\right)\\
                      \left(1-\sqrt{\mu_4},1-\sqrt{\mu_4}\right)\\
                      \left(\frac{\mu_2}{\mu_2+\mu_4}, \mu_1+\mu_3\right)\\
                      \left(\mu_1+\mu_2,\frac{\mu_3}{\mu_3+\mu_4}\right)
                      \end{tabular}
                      &
                      \begin{tabular}{L}
                      -\gamma^+ + \frac{\mu_2}{\gamma^+} + \mu_1 - \mu_4\\
                      2\sqrt{\mu_4} + \mu_1 -\mu_4 -1 \\
                      -\frac{\mu_2}{\mu_2+\mu_4}+\mu_1+\mu_2\\
                      -\frac{\mu_3}{\mu_3+\mu_4}+\mu_1+\mu_3
                      \end{tabular}& \tfour  \\\hline
5& -p-q+\mu_1-\mu_4+1 &\begin{tabular}{L}
                     (1-p) q-\mu_3\\
                     \mu_1+\mu_3-q\\
                     \mu_2-p (1-q)
                     \end{tabular}&
                     \begin{tabular}{L}
                      \left(\frac{\mu_1}{\mu_1+\mu_3},\mu_1+\mu_3\right)\\
                      \left(\frac{\mu_2}{\mu_2+\mu_4},\mu_1+\mu_3\right)
                     \end{tabular}& 
                     \begin{tabular}{L}
                      -\frac{\mu_1}{\mu_1+\mu_3} + \mu_1+\mu_2 \\
                      -\frac{\mu_2}{\mu_2+\mu_4} + \mu_1+\mu_2 
                     \end{tabular}& \tfive \\\hline
17& -p-q+\mu_1-\mu_4+1&\begin{tabular}{L}
                     \mu_1-p q\\
                     \mu_2- p (1-q)\\
                     \mu_3-(1-p) q 
                     \end{tabular}&
                     \begin{tabular}{L}
                      \left(\gamma^-,\frac{\mu_3}{1-\gamma^-}\right)\\
                      \left(\frac{\mu_1}{\mu_1+\mu_3},\mu_1+\mu_3\right)\\
                      \left(\mu_1+\mu_2,\frac{\mu_1}{\mu_1+\mu_2}\right)
                     \end{tabular}& 
                     \begin{tabular}{L}
                      -\gamma^- - \frac{\mu_3}{1-\gamma^-} +\mu_1-\mu_4+1\\
                      -\frac{\mu_1}{\mu_1+\mu_3}+\mu_1+\mu_2\\
                      -\frac{\mu_1}{\mu_1+\mu_2}+\mu_1+\mu_3
                     \end{tabular}& \tseventeen \\\hline\hline 
\multicolumn{6}{|l|}{\quad Fourth affine piece}\\\hline 
6& p+q-\mu_1+\mu_4-1 & \begin{tabular}{L}
                     \mu_2-p (1-q)\\
                     p-\mu_1-\mu_2\\
                     (1-p) q-\mu_3
                     \end{tabular}&\begin{tabular}{L}
                      \left(\mu_1+\mu_2,\frac{\mu_3}{\mu_3+\mu_4}\right)\\
                      \left(\mu_1+\mu_2,\frac{\mu_1}{\mu_1+\mu_2}\right)
                     \end{tabular}& 
                     \begin{tabular}{L}
                     \frac{\mu_3}{\mu_3+\mu_4} -\mu_1-\mu_3\\
                     \frac{\mu_1}{\mu_1+\mu_2}-\mu_1-\mu_3
                     \end{tabular} & \tsix \\\hline
7& p+q-\mu_1+\mu_4-1 & \begin{tabular}{L}
                     \mu_2-p (1-q)\\
                     \mu_3-(1-p) q\\
                     \mu_4-(1-p) (1-q)
                     \end{tabular}& \begin{tabular}{L}
                      \left(\gamma^+,\frac{\mu_3}{1-\gamma^+}\right)\\
                      \left(\frac{\mu_2}{\mu_2+\mu_4},\mu_1+\mu_3\right)\\
                     \left(\mu_1+\mu_2,\frac{\mu_3}{\mu_3+\mu_4}\right)
                     \end{tabular}& 
                     \begin{tabular}{L}
                      \gamma^++\frac{\mu_3}{1-\gamma^+}-\mu_1+\mu_4-1\\
                     \frac{\mu_2}{\mu_2+\mu_4}-\mu_1-\mu_2\\
                      \frac{\mu_3}{\mu_3+\mu_4}-\mu_1-\mu_3
                     \end{tabular} & \tseven \\\hline
8& p+q-\mu_1+\mu_4-1&\begin{tabular}{L}
                    \mu_3-(1-p) q\\
                    q-\mu_1-\mu_3\\
                    p (1-q)-\mu_2
                     \end{tabular}& \begin{tabular}{L}
                      \left(\frac{\mu_2}{\mu_2+\mu_4},\mu_1+\mu_3\right)\\
                      \left(\frac{\mu_1}{\mu_1+\mu_3},\mu_1+\mu_3\right)
                     \end{tabular}& 
                     \begin{tabular}{L}
                      \frac{\mu_2}{\mu_2+\mu_4} - \mu_1-\mu_2 \\
                      \frac{\mu_1}{\mu_1+\mu_3} - \mu_1-\mu_2 
                     \end{tabular} & \teight \\\hline
11& p+q-\mu_1+\mu_4-1&\begin{tabular}{L}
                    p q-\mu_1\\
                    p (1-q)-\mu_2\\
                    (1-p) q -\mu_3
                     \end{tabular}& \begin{tabular}{L}
                      \left(\gamma^-,\frac{\mu_3}{1-\gamma^-}\right)\\
                      \left(\sqrt{\mu_1},\sqrt{\mu_1}\right)\\
                      \left(\mu_1+\mu_2,\frac{\mu_1}{\mu_1+\mu_2}\right)\\
                      \left(\frac{\mu_1}{\mu_1+\mu_3},\mu_1+\mu_3\right)
                     \end{tabular}& 
                     \begin{tabular}{L}
                      \gamma^- + \frac{\mu_3}{1-\gamma^-} -\mu_1+\mu_4-1\\
                      2\sqrt{\mu_1} -\mu_1+\mu_4-1\\
                      \frac{\mu_1}{\mu_1 +\mu_2}-\mu_1-\mu_3\\
                      \frac{\mu_1}{\mu_1+\mu_3}-\mu_1-\mu_2
                     \end{tabular} & \televen \\\hline\hline
            \multicolumn{6}{|l|}{\quad $\gamma^+:=(1+m_2-m_3)/2 + \sqrt{(1+m_2-m_3)^2/4 - m_2}$ and $\gamma^-:=(1+m_2-m_3)/2 - \sqrt{(1+m_2-m_3)^2/4 - m_2}$ }\\\hline 
  \end{tabular}
}  
\medskip
\caption{Algebraic analysis of the Wasserstein distance function shown in Figure \ref{fig:indep_mod1}. 
}
 \label{tab:bigg1}
  \end{center}
  \end{table}                     



\subsubsection*{Acknowledgments}
GM has received funding from the European Research Council (ERC) under the European Union's Horizon 2020 research and innovation programme (grant n\textsuperscript{o}~757983).


\begin{thebibliography}{10}
\providecommand{\url}[1]{\texttt{#1}}
\providecommand{\urlprefix}{URL }
\providecommand{\doi}[1]{https://doi.org/#1}

\bibitem{MLElatentclass}
Allman, E., Cervantes, H., Evans, R., Ho\c sten, S., Kubjas, K., Lemke, D., 
Rhodes,~J., Zwiernik,~P.: {\em Maximum likelihood
  estimation of the latent class model through model boundary decomposition},
  Journal of Algebraic Statistics {\bf 34} (2019) 51--84.
  


\bibitem{arjovsky2017wasserstein}
Arjovsky, M., Chintala, S., Bottou, L.: {\em Wasserstein GAN},
{\tt arXiv:1701.07875}.

\bibitem{BASSETTI20061298}
Bassetti, F., Bodini, A., Regazzini, E.:
{\em On minimum {K}antorovich distance
  estimators}, Statistics and Probability Letters {\bf 76} (2006) 1298--1302

\bibitem{e50e1a8b234c4106af1f9dd6ef8fea08}
Bernton, E., Jacob, P., Gerber, M., Robert, C.: {\em On parameter estimation with
  the Wasserstein distance}, Information and Inference: A Journal of the IMA,
18.02.2019.

\bibitem{cuturi13}
Cuturi, M.: {\em Sinkhorn distances: Lightspeed computation of optimal transport},
 Advances in Neural Information Processing Systems,
Proceedings NIPS 2013, pp. 2292--2300.

\bibitem{DLRS}
De~Loera, J.A., Rambau, J., Santos, F.: {\em Triangulations:
Structures for Algorithms and Applications}, Algorithms and
Computation in Mathematics, vol.~25, Springer-Verlag, Berlin, 2010.

\bibitem{duarte2019discrete}
Duarte, E., Marigliano, O., Sturmfels, B.: {\em Discrete statistical models with
  rational maximum likelihood estimator},  {\tt arXiv:1903.06110}.

\bibitem{JK}
Kulas, K., Joswig, M.: {\em Tropical and ordinary convexity combined},
Advances in Geometry {\bf 10} (2010) 333--352.

\bibitem{L}
Lasserre, J.: {\em An Introduction to Polynomial and Semi-Algebraic Optimization}, Texts in Applied Mathematics, Cambridge University Press, 2015. 

\bibitem{MS}
Miller, E., Sturmfels, B.: {\em Combinatorial Commutative Algebra}, 
Graduate Texts in  Mathematics, vol.~227. Springer-Verlag, New York, 2005.

\bibitem{NR}
Nie, J., Ranestad, K.: {\em Algebraic degree of polynomial optimization},
 SIAM J.~Optim.~{\bf 20} (2009)  485--502.

\bibitem{5459199}
{Pele}, O., {Werman}, M.: {\em Fast and robust earth mover's distances}, In: 2009
  IEEE 12th International Conference on Computer Vision, pp. 460--467 (Sep 2009). 
  
\bibitem{MAL-073}
Peyre, G., Cuturi, M.:  {\em Computational optimal transport},
Foundations and Trends in Machine Learning  {\bf 11} (2019) 355--607.

\bibitem{RosStu10}
Rostalski, P., Sturmfels, B.: {\em Dualities in convex algebraic geometry}, Rendiconti di Matematica {\bf 30} (2010) 285--327.


\bibitem{seigal2018mixtures}
Seigal, A., Mont{\'u}far, G.:  {\em Mixtures and products in two graphical models},
  J.~Alg.~Stat.~{\bf 9} (2018) 1--20.

\bibitem{St_GB}
Sturmfels, B.: {\em Gr\"{o}bner Bases and Convex Polytopes}, University Lecture
  Series, vol.~8. American Mathematical Society, Providence, RI, 1996.

\bibitem{sullivant2018algebraic}
Sullivant, S.:  {\em Algebraic Statistics}, Graduate Studies in Math., American
  Mathematical Society, 2018.


\bibitem{villani08}
Villani, C.: {\em Optimal Transport: Old and New}, Grundlehren series, vol.~338, Springer 
Verlag,  2008. 

\bibitem{weed17}
Weed, J., Bach, F.: {\em Sharp asymptotic and finite-sample rates of convergence of
  empirical measures in Wasserstein distance}, {\tt  arXiv:1707.00087}.

\end{thebibliography}
\end{document}